# MACAULAY STYLE FORMULAS FOR SPARSE RESULTANTS

CARLOS D' ANDREA

Abstract. We present formulas for computing the resultant of sparse polynomials as a quotient of two determinants, the denominator being a minor of the numerator. These formulas extend the original formulation given by Macaulay for homogeneous polynomials.

## 1. Introduction

Let $\mathcal{A}_0, \ldots, \mathcal{A}_n$ be finite subsets of $\mathbb{Z}^n$ and consider $n+1$ polynomials $f_0, \ldots, f_n$ in $n$ variables such that $\mathrm{supp}(f_i) \subset \mathcal{A}_i$, $i = 0, \ldots, n$. The sparse resultant is an irreducible polynomial in the coefficients of $f_0, \ldots, f_n$, which vanishes if the system $f_i = 0, i = 0, \ldots, n$ has a solution in an algebraically closed field. It will be denoted by $\mathrm{Res}_{\mathcal{A}}(f_0, \ldots, f_n)$, where $\mathcal{A} := (\mathcal{A}_0, \ldots, \mathcal{A}_n)$.

Resultants eliminate the input variables, so they are also called eliminants. They have been used in the last decade as a computational tool for elimination of variables and for the study of complexity aspects of polynomial system solving. This has renewed the interest in finding explicit formulas for their computation (see [AS, Can1, Can2, CE1, CE2, CDS, CLO, DD, Emi1, EM, KPS, Laz, Ren, Roj, Stu1, Stu2, Stu3, ZCG]).

The study of resultants goes back to the classical work of Sylvester, Bezout, Cayley, Macaulay and Dixon in the context of homogeneous polynomials ( [Syl, Bez, Cay, Mac, Dix]). The sparse resultant, a generalization of the classical one, first appeared in the study of hypergeometric functions and $\mathcal{A}$-discriminants done by Gelfand, Kapranov and Zelevinski a few decades ago ([GKZ1, GKZ2]).

The first effective method for computing the sparse resultant was proposed by Sturmfels in [Stu1]. In [CE1, CE2], Canny and Emiris gave algorithms for computing square Sylvester style matrices with determinants equal to non-zero multiples of the resultant. By a *Sylvester*

Partially supported by Universidad de Buenos Aires, grant TX094, and Agencia Nacional de Promoción Científica y Tecnológica (Argentina), grant 3-6568.

*style matrix* we mean the matrix in the monomial bases of a linear map given by a formula as follows:

$$\begin{array}{ccc} \mathcal{S}_{\mathcal{E}_0} \oplus \mathcal{S}_{\mathcal{E}_1} \oplus \ldots \oplus \mathcal{S}_{\mathcal{E}_n} & \to & \mathcal{S}_{\mathcal{E}} \\ (g_0, g_1, \ldots, g_n) & \mapsto & \sum g_i \, f_i. \end{array}$$

Here, $\mathcal{E}_0, \ldots \mathcal{E}_n$, $\mathcal{E}$ are finite sets of monomials in a ring of Laurent polynomials $\mathbb{K}\left[x_1, x_1^{-1}, \ldots x_n, x_n^{-1}\right]$, and $\mathcal{S}_{\mathcal{B}}$ denotes the $\mathbb{K}-$vector space generated by $\mathcal{B}$. This construction was generalized by Sturmfels in [Stu2] and it was pointed out in [Emi1, CE2] that the extended formulas, when applied to the classical case, give Macaulay's original formulation (see [Mac]).

However, Macaulay succeeded in giving an explicit formula for the extraneous factor appearing in his own formulation, i.e. he showed that, in the classical case, the ratio

$$\frac{\det(\text{Sylvester matrix})}{\text{Res}_{\mathcal{A}}(f_0, \ldots, f_n)} \tag{1}$$

is a minor of the Sylvester matrix (see [Mac]). This was conjectured to happen in the sparse case, but no proof of it was available (see [CE1, CE2, CLO, Emi1, EM, GKZ2, Stu2]). For instance, in [CLO, Chapter 7], we may read:

*One of the major unsolved problems concerning sparse resultants is whether they can be represented as a quotient of two determinants. In the multipolynomial case, this is true by Theorem (4.9) of Chapter 3. Does this have a sparse analog? Nobody knows!*

In [GKZ2, Introduction], the following is written:

*Macaulay made another intriguing contribution to the theory by giving an ingenious refinement of the Cayley method* [Mac]. *It would be interesting to put his approach in the general framework of this book.*

In [Stu2, Corollary 3.1], we also find:

*It is an important open problem to find a more explicit formula for* $P_{\omega,\delta}$[1] *in the general case. Does there exist such a formula in terms of some smaller resultants?*

The main contribution of this paper is a positive answer to this question, i.e. a generalization of Macaulay's classical formulas to the sparse case by means of an explicit algorithm which produces square Sylvester style matrices. The determinant of each of these matrices is a non trivial multiple of $\text{Res}_{\mathcal{A}}(f_0, \ldots, f_n)$. Moreover, we succeed in describing the extraneous factor of our formulation (i.e. the ratio which appears in (1)) which again happens to be the determinant of a submatrix of the Sylvester matrix.

[1] $P_{\omega,\delta}$ is the extraneous factor

The paper is organized as follows: in Section 2, some notation and preliminaries are introduced. In Section 3, we explicitly construct Sylvester style matrices for generalized unmixed families of polynomials and prove that our algorithm produces formulas "à la Macaulay" for computing the sparse resultant in this case.

Section 3.2 deals with the general case, and may be regarded as an extension of the previous section. The algorithms are illustrated with examples at the end of both sections.

**Acknowledgements:** I am grateful to Alicia Dickenstein who brought my attention to this problem and to Ioannis Emiris for helpful comments. I also wish to express my deep gratitude to David Cox for for his thorough reading of preliminary drafts of this paper and very thoughtful suggestions for improvement.
This research began during the Long Semester Program in Symbolic Computation in Geometry and Analysis held at MSRI in the Fall Semester of 1998. I am grateful to the organizers for their help and support. I am especially grateful to Bernd Sturmfels for helpful conversations during those days.

## 2. Preliminaries

We review here some definitions and properties of convex polytopes and sparse resultants. More details and proofs can be found in [CLO, EM, GKZ2, Stu1, Stu2].

Let $\mathcal{A}_0, \ldots, \mathcal{A}_n$ be finite subsets of the lattice $\mathbb{Z}^n$. Set $m_i := \#(\mathcal{A}_i)$, $m := \sum_{i=0}^n m_i$ and $Q_i := conv(\mathcal{A}_i),\ i = 0, \ldots, n$. Here, $conv(\cdot)$ denote convex hull in $\mathcal{L}_{\mathbb{R}} := \mathcal{L} \otimes \mathbb{R}$, where $\mathcal{L}$ is the affine lattice generated by $\sum_{i=0}^n \mathcal{A}_i$.

For any subset $J \subset \{0, 1, \ldots, n\}$, consider the affine lattice generated by $\sum_{j \in J} \mathcal{A}_j$, and let $rk(J)$ be the rank of this lattice. For every $a \in \mathcal{A}_i$, we shall introduce a parameter $c_{i,a}$. Consider the family of *generic polynomials:*

$$f_i\left(x_1, \ldots, x_n\right) = \sum_{a \in \mathcal{A}_i} c_{i,a}\, x^a \ \ (i = 0, \ldots, n)\,. \tag{2}$$

Let $\mathbb{K}$ be an algebraically closed field. The vector of coefficients $\left(c_{i,a}\right)_{a \in \mathcal{A}_i}$ of such a family defines a point in the product of $\mathbb{K}-$projective spaces $\mathbb{P}_{\mathbb{K}}^{m_0-1} \times \ldots \times \mathbb{P}_{\mathbb{K}}^{m_n-1}$. Let $Z$ denote the subset of those families (2) which have a solution in $\left(\mathbb{K}^*\right)^n$. Here, $\mathbb{K}^*$ denotes the torus $\mathbb{K} \setminus \{0\}$. Finally, denote by $\overline{Z}$ the Zariski closure of $Z$ in $\mathbb{P}_{\mathbb{K}}^{m_0-1} \times \ldots \times \mathbb{P}_{\mathbb{K}}^{m_n-1}$.

**Theorem 2.1.** [GKZ2, Stu2] *The projective variety $\overline{Z}$ is irreducible and defined over $\mathbb{Q}$. Its codimension in $\mathbb{P}_{\mathbb{K}}^{m_0-1}\times\ldots\times\mathbb{P}_{\mathbb{K}}^{m_n-1}$ equals the maximum of $\#(I)-rk(I)$, where $I$ runs over all subsets of $\{0,1,\ldots,n\}$. The variety $\overline{Z}$ has codimension 1 if and only if there exists a unique family $\{\mathcal{A}_i\}_{i\in I}$ such that*

1. $rk(I)=\#(I)-1$,
2. $rk(J)\geq\#(J)$, *for each proper subset $J$ of $I$.*

**Definition 2.2.** [Stu2] If $I$ satisfies the conditions 1 and 2 given in the previous theorem, then the family $\{A_i\}_{i\in I}$ is said to be *essential.*

Note that if each $Q_i$ is $n-$dimensional, it is easy to check that the unique set satisfying both conditions is $I=\{0,1,\ldots,n\}$, so in this case $\overline{Z}$ has codimension 1.

The *sparse mixed resultant* $\mathrm{Res}_{\mathcal{A}}(f_0,\ldots,f_n)$ is defined as follows: if $codim(\overline{Z})=1$, then $\mathrm{Res}_{\mathcal{A}}(f_0,\ldots,f_n)$ is the unique (up to sign) irreducible polynomial in $\mathbb{Z}[c_{i,a}]$ which vanishes on $\overline{Z}$. If $codim(\overline{Z})\geq 2$, then $\mathrm{Res}_{\mathcal{A}}(f_0,\ldots,f_n)$ is defined to be the constant 1.

**Theorem 2.3.** [PS, Stu2] *If the family of supports $\{\mathcal{A}_0,\ldots,\mathcal{A}_n\}$ is essential, then for $i=0,\ldots,n$, the degree of $\mathrm{Res}_{\mathcal{A}}(f_0,\ldots,f_n)$ in the coefficients of $f_i$ is equal to the normalized mixed volume*

$$\mathcal{M}V(Q_0,\ldots,Q_{i-1},Q_{i+1},\ldots,Q_n):=$$
$$\frac{\sum_{J\subset\{0,\ldots,i-1,i+1,\ldots,n\}}(-1)^{n-\#(J)}vol\left(\sum_{j\in J}Q_j\right)}{vol(\mathcal{P})},$$

*where $vol(\cdot)$ stands for the euclidean volume in the real vector space $\mathcal{L}_{\mathbb{R}}$, and $\mathcal{P}$ is a fundamental lattice parallelotope in $\mathcal{L}$. In general, if there exists a (unique) subset $\{\mathcal{A}_i\}_{i\in i}$ which is essential, the sparse mixed resultant coincides with the resultant of the family $\{f_i:\ i\in I\}$, considered with respect to the lattice $\sum_{i\in I}\mathcal{A}_i$.*

**Example 2.4.** Set $\mathcal{A}_0=\mathcal{A}_1=\mathcal{A}_2=\{(0,0),\,(0,1),\,(1,0),\,(1,1)\}$, and consider the family

$$\begin{aligned}f_0&=c_{0,(0,0)}+c_{0,(1,0)}x_1+c_{0,(0,1)}x_2+c_{0,(1,1)}x_1\,x_2,\\ f_1&=c_{1,(0,0)}+c_{1,(1,0)}x_1+c_{1,(0,1)}x_2+c_{1,(1,1)}x_1\,x_2,\\ f_2&=c_{2,(0,0)}+c_{2,(1,0)}x_1+c_{2,(0,1)}x_2+c_{2,(1,1)}x_1\,x_2.\end{aligned}$$

Here, the Newton polytopes $Q_0$, $Q_1$, $Q_2$ are equal to the unit square $C:=[0,1]\times[0,1]$ whose vertices are precisely the points in the common support; and we have that

$$\deg_{\mathrm{coeff}_{f_i}}(\mathrm{Res}_{\mathcal{A}}(f_0,f_1,f_2))=\mathcal{M}V(C,C)=2,\ i=0,1,2.$$

A nice formula due to Dixon ([Dix]) allows us the computation of the resultant as follows:

$$\mathrm{Res}_{\mathcal{A}}(f_0,f_1,f_2)=\det\begin{pmatrix} c_{0,(0,0)} & c_{0,(1,0)} & c_{0,(0,1)} & c_{0,(1,1)} & 0 & 0\\ c_{1,(0,0)} & c_{1,(1,0)} & c_{1,(0,1)} & c_{1,(1,1)} & 0 & 0\\ c_{2,(0,0)} & c_{2,(1,0)} & c_{2,(0,1)} & c_{2,(1,1)} & 0 & 0\\ 0 & c_{0,(0,0)} & 0 & c_{0,(0,1)} & c_{0,(1,0)} & c_{0,(1,1)}\\ 0 & c_{1,(0,0)} & 0 & c_{1,(0,1)} & c_{1,(1,0)} & c_{1,(1,1)}\\ 0 & c_{2,(0,0)} & 0 & c_{2,(0,1)} & c_{2,(1,0)} & c_{2,(1,1)} \end{pmatrix}.$$

**Example 2.5.** Let

$$\begin{array}{rcl} \mathcal{A}_0 &=& \{(0,0),\,(2,2),\,(1,3)\},\\ \mathcal{A}_1 &=& \{(0,0),\,(2,0),\,(1,2)\},\\ \mathcal{A}_2 &=& \{(3,0),\,(1,1)\}. \end{array}.$$

Consider the family

$$\begin{array}{rcl} f_0 &=& \alpha_1+\alpha_2x_1^2\,x_2^2+\alpha_3x_1\,x_2^3,\\ f_1 &=& \beta_1+\beta_2x_1^2+\beta_3x_1\,x_2^2,\\ f_2 &=& \gamma_1x_1^3+\gamma_2x_1\,x_2. \end{array}$$

A straightforward computation shows that

$$\mathcal{M}V(Q_0,Q_1)=7, \mathcal{M}V(Q_0,Q_2)=7,\ \mathcal{M}V(Q_1,Q_2)=5,$$

and the sparse resultant equals

$$\begin{array}{l} \alpha_1^5\beta_3^7\gamma_1^6\gamma_2+3\alpha_1^4\alpha_2\beta_2^2\beta_3^5\gamma_1^4\gamma_2^3+3\alpha_1^3\alpha_2^2\beta_2^4\beta_3^3\gamma_z^2\gamma_2^5\\ -13\alpha_1^3\alpha_2\alpha_3\beta_1^2\beta_2\beta_3^4\gamma_1^5\gamma_2^2-7\alpha_1^3\alpha_3^2\beta_1\beta_2^3\beta_3^3\gamma_1^4\gamma_2^3+6\alpha_1^2\alpha_2^3\beta_1^3\beta_2\beta_3^3\gamma_1^4\gamma_2^3\\ +\alpha_1^2\alpha_2^3\beta_2^6\beta_3\gamma_2^7-\alpha_1^2\alpha_2^2\alpha_3\beta_1^2\beta_2^3\gamma_1^3\gamma_2^4+5\alpha_1^2\alpha_2\alpha_3^2\beta_1^4\beta_3^3\gamma_1^6\gamma_2\\ -\alpha_1^2\alpha_2\alpha_3^2\beta_1\beta_2^5\beta_3\gamma_1^2\gamma_2^5+14\alpha_1^2\alpha_3^3\beta_1^3\beta_2^2\beta_3^2\gamma_1^5\gamma_2^2+\alpha_1^2\alpha_3^3\beta_2^7\gamma_1\gamma_2^6\\ -2\alpha_1\alpha_2^4\beta_1^3\beta_2^3\beta_3\gamma_1^2\gamma_2^5-5\alpha_1\alpha_2^3\alpha_3\beta_1^5\beta_3^2\gamma_1^5\gamma_2^2+\alpha_2^5\beta_1^6\beta_3\gamma_1^4\gamma_2^3\\ +2\alpha_1\alpha_2^2\alpha_3^2\beta_1^4\beta_2^2\beta_3\gamma_1^4\gamma_2^3-2\alpha_1\alpha_2\alpha_3^3\beta_1^3\beta_2^4\gamma_1^3\gamma_2^4-7\alpha_1\alpha_3^4\beta_1^5\beta_2\beta_3\gamma_1^6\gamma_2\\ +\alpha_2^2\alpha_3^3\beta_1^6\beta_2\gamma_1^5\gamma_2^2+\alpha_3^5\beta_1^7\gamma_1^7 \end{array}$$

For an explicit computation of this resultant, see [Stu2].

As usual, to define a *face* of a polytope $Q\subset\mathcal{L}_{\mathbb{R}}\subset\mathbb{R}^n$, let $v$ be a vector in $\mathbb{R}^n$. Set

$$m_Q(v):=\min_{q\in Q}\{\langle q,v\rangle\},$$

and call

$$Q_v:=Q\cap\{m\in\mathbb{R}^n:\ \langle m,v\rangle=m_Q(v)\}$$

*the face of* $Q$ *determined by* $v$. The vector $v$ will be called an *inward normal vector* of $Q_v$. If $\dim(Q_v)=n-1$, then $Q_v$ will be called *facet.*

## 3. Macaulay style formulas for generalized unmixed families of polynomials

In this section, we will construct square Sylvester style matrices in the case all $Q_i$ are integer multiples of a fixed polytope $P$, i.e. there exists positive integer numbers $k_0, k_1, \ldots, k_n$ such that

$$Q_i = conv(\mathcal{A}_i) = k_i\, P,\ i = 0, \ldots, n.$$

This case is treated in [CDS] in a toric setting. These families may be identified with homogeneous polynomials in the coordinate ring of a projective toric variety (in the sense of [Cox]), with "degrees" $\alpha_0, \ldots, \alpha_n$, where every $\alpha_i$ is $\mathbb{Q}-$ample (see [Ful]). We shall call them *generalized unmixed polynomials* because they contain the well-known unmixed family of polynomials, which is the case when all input supports $\mathcal{A}_i$ coincide.

The matrices to be constructed here also generalize the formulas given by Macaulay in [Mac] in the homogeneous case, where all the supports are multiples of the standard simplex

$$S_n := \{(q_1, \ldots, q_n) \in \mathbb{R}^n :\ 0 \leq q_i \leq 1,\ \sum_{i=1}^{n} q_i \leq 1\}.$$

*Warning*: Some care must be taken with taking convex hulls, because the sparse resultant depends strongly on the finite data $(\mathcal{A}_0, \mathcal{A}_1, \ldots, \mathcal{A}_n)$ and different families of input supports may give the same polytopes $Q_0, \ldots, Q_n$ (see [Stu2]).

*Remark* 3.1. It is easy to see that, in the generalized unmixed case, because of Theorem 2.1, in order to have a nontrivial resultant, $P$ must be $n-$dimensional. It is also clear that, in this case, $\mathcal{L}_{\mathbb{R}} = \mathbb{R}^n$.

Given $\lambda \in \mathbb{Q}_{\geq 0}$ and a generic $\delta \in \mathcal{L} \otimes \mathbb{Q}$ as in [CE1, CE2, Stu2], our algorithm will produce a Sylvester style matrix $\mathbb{M}$ whose rows and columns will be indexed by the integer points in

$$\mathcal{E} := ((k_0 + k_1 + \ldots + k_n + \lambda)P + \delta) \cap \mathcal{L}, \tag{3}$$

and whose determinant will be a nonzero multiple of the sparse resultant. Moreover, we shall be able to identify the extraneous factor $\det(\mathbb{M})/\mathrm{Res}_{\mathcal{A}}(f_0, \ldots, f_n)$ as a minor of this determinant.

The algorithm is recursive in the dimension of the polytope $P$, and in its intermediate steps, uses the mixed-subdivision algorithm of Canny and Emiris (see, [CE1, CE2, Stu2]) in order to refine the subdivision (see the comments in Remark 3.15).

*Remark* 3.2. It was stated in the introduction that $\mathcal{E}$ should be a set of monomials. Indeed this is true provided that we identify an integer point $\alpha \in \mathcal{E}$ with the Laurent monomial $x^\alpha$.

3.1. **Constructing the Matrix $\mathbb{M}$.** Given $\lambda$ and $\delta$ as before, set $Q := \lambda\, P$. Let $V(Q) \subset \mathcal{L} \otimes \mathbb{Q}$ be the set of vertices of $Q$. Observe that they are not necessarily integer points.

Choose a vertex $\mathbf{b_0} \in \mathcal{A}_0$, and consider the following lifting functions:

$$\begin{array}{rccll} \omega_0: & \mathcal{A}_0 & \to & \mathbb{R} & \\ & \mathbf{b_0} & \mapsto & 1 & \\ & b & \mapsto & 0 & \text{if}\, b \neq \mathbf{b_0} \\ \\ \omega_i: & \mathcal{A}_i & \to & \mathbb{R} & \\ & b & \mapsto & 0 & \forall b \in \mathcal{A}_i, \quad i = 1, \ldots, n \\ \\ \omega: & V(Q) & \to & \mathbb{R} & \\ & b & \mapsto & 0 & \forall b \in V(Q). \end{array} \tag{4}$$

Set

$$\Omega := \left(\omega_i(b)_{b \in \mathcal{A}_i,\, i=0,1,\ldots,n},\ \omega(b)_{b \in V(Q)}\ \right) \in \mathbb{R}^{m+\#V(Q)}$$

and consider the lifted polytopes in $\mathbb{R}^{n+1}$ :

$$\begin{array}{rcl} Q_{i,\Omega} & := & conv\{(a, \omega_i(a)) :\, a \in \mathcal{A}_i\} \\ Q_\Omega & := & conv\{(b, \omega(b)) :\, b \in V(Q)\}. \end{array}$$

By projecting the upper envelope of $Q_{i,\Omega}$ (resp. $Q_\Omega$) we get a *coherent mixed decomposition* $\Delta_{i,\Omega}$ (resp. $\Delta_\Omega$) of the polytopes $Q_i$ (resp. $Q$). The cells in this decomposition are the projections of precisely those faces of $Q_{i,\Omega}$ and $Q_\Omega$ on which a linear functional with negative last coordinate is minimized (see [Stu2]).

Similarly, by projecting the upper envelope of the Minkowski sum

$$Q_{0,\Omega} + Q_{1,\Omega} + \ldots + Q_{n,\Omega} + Q_\Omega,$$

we get a coherent mixed decomposition of

$$\mathcal{Q} := Q_0 + Q_1 + \ldots + Q_n + Q. \tag{5}$$

Each cell in this decomposition is of the form

$$F = F_0 + F_1 + \ldots + F_n + F,$$

where $F_i$ (resp. $F$) is a cell in $\Delta_{i,\Omega}$, (resp. $\Delta_\Omega$).

Because the lifting functions we have used in (4) are most of them trivial, we can characterize all the proper cells in $\Delta_{1,\Omega}, \ldots, \Delta_{n,\Omega}, \Delta_\Omega$ :

- the whole polytopes $Q_1, \ldots, Q_n, Q$ corresponding to the linear functional associated to the vector $(0, -1) \in \mathbb{R}^n \times \mathbb{R}$;

- for every $v \in \mathbb{R}^n \setminus \{0\}$, the faces $Q_{1v}, \dots, Q_{nv}, Q_v$ associated to the vector $(v, \alpha)$, where $\alpha$ is any negative number.

On the other hand, on $\Delta_{0,\Omega}$, the cell corresponding to $(0, -1)$ is the singleton $\{\mathbf{b_0}\}$, and it is easy to check that every cell of dimension $n$ in its decomposition is the convex hull of $\mathbf{b_0}$ and a facet of $Q_0$ not containing this point. We shall call it $F_{0,v}$, where $v$ denotes the integer primitive inward normal of the facet.

So, we can characterize all maximal (i.e. of dimension $n$) cells in the polyhedral decomposition of $\mathcal{Q}$ as follows:

- $\{\mathbf{b_0}\} + Q_1 + \ldots + Q_n + Q$. This shall be called the *primary* cell.
- $F_{0,v} + Q_{1v} + \ldots + Q_{nv} + Q_v$, for some $v \in \mathbb{R}^n$. These will be called *secondary* cells, and are associated to a nonzero vector $v \in \mathbb{R}^n$.

The following Lemma will be useful in the sequel:

**Lemma 3.3.** *Suppose that* $\dim(F_{0,v}) = n$. *Then, for every* $q \in \mathbb{Q}$, *if* $F_{0,v} \cap \{m \in \mathcal{L}_{\mathbb{R}} : \langle m, v\rangle = q\} \neq \emptyset$, *there is* $\lambda_q \in \mathbb{Q}_{\geq 0}$ *such that the intersection is a polytope congruent to* $\lambda_q P_v$.

*Proof.* Suppose w.l.o.g. that $\mathbf{b_0} = 0$, then $F_{0,v}$ is the convex hull of the origin in $\mathcal{L}_{\mathbb{R}}$ and a finite set of points $\{v_1, \dots, v_M\}$ all of them satisfying $\langle v_i, v\rangle = \lambda_v < 0$, hence lying in a hyperplane $\mathcal{H}_v$ not passing through the origin.

The intersection of $F_{0,v}$ with a hyperplane parallel to $\mathcal{H}_v$ will be nonempty if and only if $\lambda_v \leq q \leq 0$. If this happens, the intersection will be the convex hull of

$$\{\frac{q}{\lambda_v} v_1, \dots, \frac{q}{\lambda_v} v_M\},$$

which is equal to $\frac{q}{\lambda_v} Q_{0v} = \frac{q}{\lambda_v} k_0 P_v$. □

In order to construct the Sylvester style matrix, we will take into account whether the points lie in a translation of the primary cell or not. The first requirement we impose on $\delta$ is that every point in $(\mathcal{Q} + \delta) \cap \mathcal{L}$ must belong to the interior of a shifted maximal cell (primary or secondary).

3.1.1. *Points in the shifted primary cell.* Proceed as in [CE1, CE2, Emi1, Stu2]: choose generic lifting functions $\tilde{\omega}_1, \dots, \tilde{\omega}_n, \tilde{\omega}$, defined over $\mathcal{A}_1, \dots, \mathcal{A}_n$, $V(Q)$, in such a way that they produce a *tight mixed coherent decomposition* of the Minkowski sum

$$Q_1 + \ldots + Q_n + Q.$$

This implies that each $n-$ dimensional cell $\mathcal{F}$ in this decomposition equals

$$F_1 + F_2 + \ldots + F_n + F,$$

where $F_i$ is a cell in $\Delta_{i,\tilde{\omega}}$, $F$ is a cell in $\Delta_{\tilde{\omega}}$ and

$$n = \dim(F_1) + \ldots + \dim(F_n) + \dim(F).$$

As a consequence, at least one of these dimensions is equal to 0. The *row content* of $p \in (\{\mathbf{b_0}\} + Q_1 + \ldots + Q_n + Q + \delta) \cap \mathcal{L}$ is a pair $(i, a)$ defined as follows:

- if $p - \delta - \mathbf{b_0}$ lies in the cell $\mathcal{F} = F_1 + \ldots + F_n + F$, and the set $\{j : 1 \leq j \leq n, \dim(F_j) = 0\}$ is not empty, let $i$ be the largest index such that $\dim(F_i) = 0$, and let $F_i = \{a\}$.
- If $\dim(F_j) > 0$, for all $j$, (because of the genericity of the lifting functions $\tilde{\omega}_j$, this implies that $\dim(F) = 0$) then $i := 0$, and $a := \mathbf{b_0}$.

We shall say that $\mathcal{F}$ is *mixed of type* 0 if the last item holds; otherwise, the cell shall be called *non-mixed.*

*Remark* 3.4. The concepts of *row content, mixed and non-mixed cells* defined previously, appear with a slightly different meaning in [CE1, CE2, Stu2]. We shall discuss some relations between both definitions in example 3.3.3.

We can now fill the rows of the matrix $\mathbb{M}$ indexed by those points $p$ lying in $(\{\mathbf{b_0}\} + Q_1 + \ldots + Q_n + Q + \delta) \cap \mathcal{L}$ as follows: for every $p' \in \mathcal{E}$, the entry indexed by $(p, p')$ equals the coefficient of $x^{p'}$ in the expansion of the polynomial $x^{p-a} f_i(x)$. Here, $(i, a)$ is the row content of $p$.

3.1.2. *Points in the shifted secondary cells.* Let $v \in \mathbb{R}^n$ such that

$$\mathcal{F}_v = F_{0,v} + Q_{1_v} + \ldots + Q_{n_v} + Q_v \tag{6}$$

is a maximal cell in the polyhedral decomposition of $\mathcal{Q}$. This implies that

- $\dim(P_v) = n - 1$, $Q_v = \lambda P_v$, and $Q_{i_v} = k_i P_v$, $i = 0, 1, \ldots, n$.
- $\dim(F_{0,v}) = n$.

Intersecting the shifted cell $F_{0,v}$ with hyperplanes parallel to $P_v$ and using Lemma 3.3, it is easy to see that $(\mathcal{F}_v + \delta) \cap \mathcal{L}$ can be written as a disjoint union of sets of the type

$$\Big((\tilde{\lambda} P_v + \delta' + k_1 P_v + \ldots + k_n P_v + \lambda P_v + \delta\Big) \cap \mathcal{L},\ \tilde{\lambda} \in \mathbb{Q}_{\geq 0},$$

which may be rearranged as follows:

$$(k_1 P_v + \ldots + k_n P_v + \lambda_v P_v + \delta_\lambda) \cap \mathcal{L},\ \lambda_v \in \mathbb{Q}_{\geq 0}, \tag{7}$$

where $\delta_\lambda := \delta + \delta'$. Here is where the recursion step comes: consider the $v-$*facet family*

$$f_{i_v}(x_1, \ldots, x_n) = \sum_{a \in Q_{i_v} \cap \mathcal{A}_i} c_{i,a}\, x^a \ (i = 1, \ldots, n). \tag{8}$$

Due to the fact that $Q_{i_v} = k_i\, P_v$, has dimension exactly $n-1$, it is straightforward to check that the family $\{Q_{i_v} \cap \mathcal{A}_i\}_{i=1,\ldots,n}$ is essential. This implies that the sparse resultant of the polynomials (8) is not constantly equal to one. We shall denote this sparse resultant as

$$\mathrm{Res}_v(f_{1_v}, \ldots, f_{n_v}). \tag{9}$$

In order to use the inductive hypothesis, we must decrease the dimension of the supports with some care: let $\mathcal{L}_v \subset \mathcal{L}$ be the lattice which is orthogonal to $v$, and denote by $\mathcal{L}_{\mathcal{A}_{1_v}+\cdots+\mathcal{A}_{n_v}}$ the affine lattice generated by $\mathcal{A}_{1_v} + \cdots + \mathcal{A}_{n_v}$.

After a translation, we may suppose w.l.o.g. that $0 \in \mathbb{R}^n$ is a vertex of $P_v$. This implies that $\mathcal{L}_{\mathcal{A}_{1_v}+\cdots+\mathcal{A}_{n_v}}$ is a sublattice of $\mathcal{L}_v$, both having dimension $n-1$, and we may consider the index $[\mathcal{L}_v : \mathcal{L}_{\mathcal{A}_{1_v}+\cdots+\mathcal{A}_{n_v}}]$, which will be denoted by $\mathrm{ind}_v$. Let $\mathrm{q}_1, \cdots, \mathrm{q}_{\mathrm{ind}_v}$ be coset representatives for $\mathcal{L}_{\mathcal{A}_{1_v}+\cdots+\mathcal{A}_{n_v}}$ in $\mathcal{L}_v$.

For every $p \in \mathcal{L}$, there exists a unique $j \in \{1, \ldots, \mathrm{ind}_v\}$ such that $p \in (\mathrm{q}_j + \mathcal{L}_{\mathcal{A}_{1_v}+\cdots+\mathcal{A}_{n_v}}) \oplus \mathbb{Z}v$, so $p$ may be written as $p_v + p^v$, the latter being an integer multiple of $v$, and $p_v \in \mathrm{q}_j + \mathcal{L}_{\mathcal{A}_{1_v}+\cdots+\mathcal{A}_{n_v}}$. Also, $\delta_\lambda \in \mathcal{L} \otimes \mathbb{Q} = (\mathcal{L}_{\mathcal{A}_{1_v}+\cdots+\mathcal{A}_{n_v}} \oplus \mathbb{Z}v) \otimes \mathbb{Q}$ may be decomposed as $\delta_\lambda^v + \delta_{\lambda v}$, where $\delta_\lambda^v$ (resp. $\delta_{\lambda v}$) lies in $\mathbb{Q}v$ (resp. $\mathcal{L}_{\mathcal{A}_{1_v}+\cdots+\mathcal{A}_{n_v}} \otimes \mathbb{Q}$).

If, in addition, $p$ belongs to (7), we must have $p^v = \delta_\lambda^v$. This is due to the fact that $k_1 P_v + \ldots + k_n P_v + \lambda_v P_v + \delta_{\lambda v} \subset \mathcal{L}_{\mathcal{A}_{1_v}+\cdots+\mathcal{A}_{n_v}} \otimes \mathbb{R}$. So, $p - p^v = p - \delta_{\lambda v}$ lies in

$$(k_1 P_v + \ldots + k_n P_v + \lambda_v P_v + \delta_{\lambda v}) \cap (\mathrm{q}_j + \mathcal{L}_{\mathcal{A}_{1_v}+\cdots+\mathcal{A}_{n_v}}).$$

Finally, set $\delta_{jv} := \delta_{\lambda v} - \mathrm{q}_j \in \mathcal{L}_{\mathcal{A}_{1_v}+\cdots+\mathcal{A}_{n_v}} \otimes \mathbb{Q}$. Now it is straightforward to check that a point $p$ belongs to (7), if and only if, there exists $j = 1, \ldots, \mathrm{ind}_v$, such that $\mathbf{p} := p - \delta_\lambda^v - \mathrm{q}_j$ lies in

$$(k_1 P_v + \ldots + k_n P_v + \lambda_v P_v + \delta_{jv}) \cap \mathcal{L}_{\mathcal{A}_{1_v}+\cdots+\mathcal{A}_{n_v}}. \tag{10}$$

We have decreased dimension, and may use the inductive hypothesis in order to compute the resultant (9) using $\mathcal{L}_{\mathcal{A}_{1_v}+\cdots+\mathcal{A}_{n_v}}$ instead of $\mathcal{L}$ and $\delta_{jv}$ instead of $\delta$. It turns out that, for every $\lambda_v \geq 0$, there will be $\mathrm{ind}_v$ square matrices of the form $\mathbb{M}_{v,\lambda_v}$ indexed by the points $\mathbf{p}$ lying in (10). Using the monomial bijection $p = \mathbf{p} + \delta_\lambda^v + \mathrm{q}_j$, we can relabel rows and columns of these matrices with the points of (7). The determinant of each $\mathbb{M}_{v,\lambda_v}$ will be a non trivial multiple of (9), and each of these determinants will have the same degree as $\mathrm{Res}_v(f_{1_v}, \ldots, f_{n_v})$ in the coefficients of $f_{1_v}$ ($f_{1_v}$ will play here the role of $f_0$ in the previous step).

Again, as in the primary cell, each of these matrices has in the row corresponding to a point $p$, the coordinates of the monomial expansion of $x^{p-a}\, f_{i_v}(x)$ for some $a$, $i$. In order to fill the row of $\mathbb{M}$ corresponding to $p$, we proceed as before: the entry indexed by $(p, p')$ equals the coefficient of $x^{p'}$ in the expansion of the polynomial $x^{p-a}\, f_i(x)$.

In order to finish the algorithm properly, the reader should check that, in the case $n = 1$, this procedure constructs a classical "Sylvester style matrix" ([Syl, Mac], see also Example 3.3.1) for two polynomials in one variable. If the supports generate the affine lattice $\mathbb{Z}$, the matrix will be indexed by a set of monomials of the type

$$\{x_1^a,\, x_1^{a+1}, \ldots, x_1^{a+s}\}\ a \in \mathbb{Z}, s \in \mathbb{N}.$$

*Remark* 3.5. Observe that, at each step of the recursion, we need to impose some conditions on the different $\delta$'s in order to guarantee that all integer points are in the interior of a cell in each intermediate step. This happens for $\delta$ generic.

*Remark* 3.6. Let $d_v$ be the *$v-$lattice diameter* of the cell $\mathcal{F}_v$ which is defined as follows:

$$d_v := \max_{m\in\mathcal{F}_v} \langle m, v\rangle - \min_{m\in\mathcal{F}_v} \langle m, v\rangle. \tag{11}$$

It is straightforward to check that, the number of matrices of the type $\mathbb{M}_{v,\lambda_v}$ is exactly $\mathbf{d_v} := d_v\, \mathrm{ind}_v$. Besides, due to the fact that $F_1, F_2, \ldots, F_n, F$ are facets associated to $v$, the difference (11) may be actually computed as

$$\max_{m\in F_{0,v}} \langle m, v\rangle - \min_{m\in F_{0,v}} \langle m, v\rangle. \tag{12}$$

*Remark* 3.7. At each step of the recursion, the shifted secondary cells of the previous step are partitioned in such a way that their integer points are distributed into new primary and secondary cells. The new primary is again subdivided into mixed and unmixed cells. We shall keep track of this information, so a point will be said to be in a mixed cell of type $i$, if it belongs to a mixed cell which appeared at step $i+1$. It is easy to see that, if a point is in a mixed cell of type $i$, then the row it indicates contains some coefficients of the expansion of a multiple of $f_i$.

It is also clear that, at the end of the recursion, each point in $\mathcal{E}$ has associated a row content.

3.2. **Generalized Macaulay Formula.** Now we are ready to state and prove the central result of this section. Before doing that, note that $\det(\mathbb{M})$ is well defined up to sign, because we have not given any order among the elements of $\mathcal{E}$. Also, $\mathrm{Res}_{\mathcal{A}}(f_0, \ldots, f_n)$ is well defined up to sign. So, the following statement will be true up to sign.

**Theorem 3.8.** $\mathbb{M}$ *is a generically nonsingular Sylvester style matrix. Moreover, we have the following formula "à la Macaulay":*

$$\det(\mathbb{M}) = \mathrm{Res}_{\mathcal{A}}(f_0, \ldots, f_n)\, \det(\mathbb{E}),$$

*where* $\mathbb{E}$ *is the square submatrix of* $\mathbb{M}$ *formed by omiting all rows and columns indexed by points lying in mixed cells.*

*Remark* 3.9. It is easy to see that $\mathbb{E}$ does not contain coefficients of $f_0$, so as an inmediate corollary of Theorem 3.8 we get that $\det(\mathbb{M})$ has the same degree as $\mathrm{Res}_{\mathcal{A}}(f_0, \ldots, f_n)$ in the coefficients of $f_0$. Replacing the role played by $f_0$ with any $f_i,\ i = 1, \ldots, n$, we also have a formula for computing $\mathrm{Res}_{\mathcal{A}}(f_0, \ldots, f_n)$ as the gcd of $n+1$ determinants (see also [CE1, CE2, EM]).

*Proof.* First of all, we will prove that, for a given point $p \in \mathcal{E}$, if $p'$ belongs to the support of $x^{p-a} f_i$, then $p'$ must also be a point of $\mathcal{E}$. Here, the pair $(i, a)$ is the row content of the point $p$. This will imply that $\mathbb{M}$ is a Sylvester style matrix.

Two different scenarios must be considered:

- If the point belongs to the shifted primary cell, proceeding as in [CE1, CE2], it is easy to see that this happens.
- If the point belongs to a shifted secondary cell, let us say $\mathcal{F}_v$, because of (6), $x^{p-a} f_{i_v}$ has its support contained in

$$F_{0,v} + Q_{1_v} + \ldots + Q_{n_v} + Q_v + \delta.$$

  This, combined with the fact that, in secondary cells, $i$ is always bigger than 0, implies that the support of $x^{p-a} f_i$ is contained in

$$F_{0,v} + Q_{1_v} + \ldots + Q_{i-1_v} + Q_i + Q_{i+1_v} + \ldots Q_{n_v} + Q_v + \delta,$$

  the latter set being a subset of $\mathcal{Q} + \delta$. From here, the claim follows straightforwardly.

Once we know that $\mathbb{M}$ is a Sylvester style matrix, it is easy to see that $\mathrm{Res}_{\mathcal{A}}(f_0, \ldots, f_n)$ divides $\det(\mathbb{M})$ using the standard argumentation given in [CE1, CE2, Stu2].

We shall regard $\det(\mathbb{M})$ as a polynomial in $\mathbb{Z}[c_{i,a}]$ and will prove that it is not identically zero by showing that its highest term with respect to some monomial order is nonzero. Explicitly, we shall prove that, for the vector

$$\overline{\omega} := (\omega_0, \omega_1, \ldots, \omega_n) \in \mathbb{R}^m, \tag{13}$$

where the $\omega_i$ were defined in (4),

$$init_{\overline{\omega}}\,(\det(\mathbb{M})) \neq 0. \tag{14}$$

It is easy to see that, writing $\det(\mathbb{M})$ as a polynomial in $c_{0,\mathbf{b_0}}$ with coefficients in $\mathbb{Z}[c_{i,a} \setminus \{c_{0,\mathbf{b_0}}\}]$, the leading term of this polynomial is (14).

Due to the special role that $f_0$ has played in the construction of the matrix, the number of integer points lying in shifted $0-$mixed cells is equal to $M_0 := \mathcal{M}V(Q_1, \ldots, Q_n)$, i.e. the degree of $\mathrm{Res}_{\mathcal{A}}(f_0, \ldots, f_n)$ in the coefficients of $f_0$ ([HS, CE1, CE2, Stu2]). So, it is straightforward to check that

$$\deg_{\mathrm{coeff}(f_0)}(\mathrm{Res}_{\mathcal{A}}(f_0, \ldots, f_n)) = \deg_{\mathrm{coeff}(f_0)}(\det(\mathbb{M}))\,.$$

This, combined with the special way in which we have lifted the polytopes (i.e. just lifting the point $\mathbf{b_0}$), implies that

$$\frac{\det(\mathbb{M})}{\mathrm{Res}_{\mathcal{A}}(f_0, \ldots, f_n)} = \frac{init_{\overline{\omega}}\,(\det(\mathbb{M}))}{init_{\overline{\omega}}\,(\mathrm{Res}_{\mathcal{A}}(f_0, \ldots, f_n))}$$

$$= \frac{\text{coefficient of}\, c_{0,\mathbf{b_0}}^{M_0} \text{ in } \det(\mathbb{M})}{\text{coefficient of}\, c_{0,\mathbf{b_0}}^{M_0} \text{ in } \mathrm{Res}_{\mathcal{A}}(f_0, \ldots, f_n)}.$$

In order to prove the theorem, we shall proceed as in [Mac], by showing that the numerator of this fraction is nonzero, and that the extraneous factor, i.e. the ratio, is $\det(\mathbb{E})$. The proof will be again by recurrence on $n$.

The basic case ($n = 1$) is completely contained in the classical formulas given by Macaulay in [Mac]. In the general case, in order to compute $init_{\overline{\omega}}\,(\det(\mathbb{M}))$, we proceed as in [CE1, CE2, Stu2]: replace the polynomials (2) by the following deformed family:

$$f_{i,\omega} := \sum_{a \in \mathcal{A}_i} c_{i,a}\, t^{\omega_i(a)} x^a, \; i = 0, 1, \ldots, n, \tag{15}$$

and consider the deformed matrix $\mathbb{M}\left(c_{i,a}\, t^{\omega_i(a)}\right)$.

*Remark* 3.10. Actually,

$$\begin{array}{rcll} f_{0,\omega} & = & c_{0,\mathbf{b_0}}\, t\, x^{\mathbf{b_0}} + \sum_{a \in \mathcal{A}_i \setminus \{\mathbf{b_0}\}} c_{i,a} x^a, & \\ f_{i,\omega} & = & f_i & i \geq 1. \end{array}$$

It is easy to see that $init_{\overline{\omega}}\,(\det(\mathbb{M}))$ is the leading coefficient of the determinant of the matrix $\mathbb{M}(c_{i,a}\, t^{\omega_i(a)})$ regarded as a polynomial in $t$.

For every $p \in \mathcal{E}$, let $h(p)$ be the biggest rational number such that

$$(p-\delta, h(p)) \in \mathcal{Q}_\Omega = Q_{0,\Omega} + Q_{1,\Omega} + \ldots + Q_{n,\Omega} + Q_\Omega.$$

The following observations will be useful later:

**Lemma 3.11.** *The function $h$ verifies*

1. $0 < h(p) \leq 1$, *for all* $p \in \mathcal{E}$.
2. $h(p) = 1$ *if and only if* $p - \delta$ *lies in the primary cell.*
3. *If* $p-\delta$ *and* $q-\delta$ *are in the same secondary cell, say* $\mathcal{F}_v$, *then*
$$h(p) = h(q) \iff \langle p, v\rangle = \langle q, v\rangle.$$
4. *If* $p \in (\mathcal{F}_v + \delta) \cap \mathcal{L}$, *the row content of* $p$ *is the pair* $(i,a)$, *and* $v' \neq \mu\, v,\ \mu > 0$, *then*
$$((p-\delta-a, h(p)) + Q_{i,\Omega}) \cap \mathcal{Q}_{\Omega(v',l')} = \emptyset,\ \forall l' \in \mathbb{R}_{<0};$$
*here,* $\mathcal{Q}_{\Omega(v',l')} \subset \mathbb{R}^{n+1}$ *is the face of* $\mathcal{Q}_\Omega$ *determined by* $(v', l')$.

*Proof of the Lemma:* The first two statements are obvious. The third assertion holds straightforwardly just noting that $\mathcal{F}_v$ is the projection of $\mathcal{Q}_{\Omega(v,l)}$, where $l$ is a negative real number, and

$$p-\delta,\ q-\delta \in \mathcal{F}_v \iff (p-\delta, h(p))\,, (q-\delta, h(q)) \in \mathcal{Q}_{\Omega(v,l)},$$

hence

$$\langle (p-\delta, h(p)), (v,l)\rangle = \langle (q-\delta, h(q)), (v,l)\rangle.$$

Lastly, let $l \in \mathbb{R}_{<0}$ such that $\mathcal{Q}_{\Omega(v,l)}$ projects bijectively onto $\mathcal{F}_v$. Due to the generic conditions imposed on $\delta$, it turns out that the point $p_\Omega := (p-\delta, h(p))$ belongs to the relative interior of $\mathcal{Q}_{\Omega(v,l)}$. If $v' = \mu\, v$, and $l' = \mu\, l$, with $\mu \leq 0$, this would imply $l' \in \mathbb{R}_{\geq 0}$, which is not of interest for us, so we can suppose w.l.o.g. that $(v', l')$ is not parallel to $(v,l)$. Then, one can slightly displace the point $p_\Omega$ inside $\mathcal{Q}_{\Omega(v,l)}$, in the direction of the orthogonal projection of $-(v', l')$ over the hyperplane $\{\langle x, (v,l)\rangle = 0\} \subset \mathbb{R}^{n+1}$. After this displacement, all the points in the shifted $Q_i$ will still lie in $\mathbb{Q}_\Omega$. This is due to the convexity argument given in [CE1, CE2], which states that for *every* point $q_\Omega$ lying in $\mathcal{Q}_{\Omega(v,l)}$,

$$q_\Omega - (a,0) + Q_{i,\Omega} \subset \mathcal{Q}_\Omega.$$

So, if the statement of the Lemma were not true, the common points in the intersection with $\mathcal{Q}_{\Omega(v',l')}$ would not belong to $\mathcal{Q}_\Omega$ after the displacement, which is impossible.

□

We shall use the convexity argument given in [CE1, CE2, Stu2] as follows: for each $p \in \mathcal{E}$, we multiply every element in the row indexed by $p$ by $t^{h(p)-\omega_i(a)}$, where the row content of $p$ is $(i,a)$. Let us call the matrix obtained in this way $\mathbb{M}'(t)$. It is easy to see that the leading coefficient of $\det(\mathbb{M}'(t))$ (as a polynomial in $t$) is $init_{\overline{\omega}}(\det(\mathbb{M}))$.

**Proposition 3.12.** *Let $0 < \gamma_1 < \gamma_2 < \ldots < \gamma_N = 1$ be the different values for $h(p)$ as $p$ ranges in $\mathcal{E}$. Then, the leading coefficient of $\det(\mathbb{M}'(t))$ (as a polynomial in $t$) factorizes as*

$$\prod_{j=1}^{N} \det\left(\mathbb{M}_j\right), \tag{16}$$

*where $\mathbb{M}_j$ is the square submatrix of $\mathbb{M}$ made by choosing all rows and columns indexed by points $p$ such that $h(p) = \gamma_j$.*

*Moreover, the product (16) can be reorganised as follows:*

$$\det(\mathbb{M}_N) \prod_{\mathcal{F}_v} \left(\det(\mathbb{M}_v^1) \ldots \det(\mathbb{M}_v^{\mathbf{d_v}})\right), \tag{17}$$

*where the product is taken over all secondary cells $\mathcal{F}_v$ and $\mathbb{M}_v^i$ is one of the matrices $\mathbb{M}_{v,\lambda_v}$ defined in Section 3.1.2, i.e. a Sylvester style matrix for the sparse resultant associated to the family (8).*

*Proof of the Proposition:* We will compute the leading term of $\det(\mathbb{M}'(t))$ by searching, in each column, the highest power of $t$ appearing in that column, replacing by 0 all the entries which have not this highest power, and computing the determinant of the modified matrix. We shall call this matrix $mod(\mathbb{M}'(t))$.

It is straightforward to check that, in $\mathbb{M}'(t)$, the highest power appearing in the column indexed by $p'$ is exactly $h(p')$. So, in $mod(\mathbb{M}'(t))$, in the columns indexed by those $p'$ such that $h(p') = \gamma_1$ there cannot be nonzero entries in the rows indexed by those $p$ such that $h(p) \neq \gamma_1$; otherwise $h(p) > \gamma_1$ and this will imply that the power of $t$ appearing in the $(p, p')$ place of the matrix will be strictly bigger than $\gamma_1$, which is impossible.

This implies that $mod(\mathbb{M}'(t))$, after division by some power of $t$, and ordering its rows and columns by putting the points $p'$ such that $h(p') = \gamma_1$ at the beginning, has the following structure:

$$\begin{pmatrix} \mathbb{M}_1 & A \\ 0 & B \end{pmatrix}. \tag{18}$$

Repeating the argument recursively, it turns out that the structure of $mod(\mathbb{M}'(t))$ is triangular as follows:

$$\begin{pmatrix} \mathbb{M}_1 & * & \dots & * \\ 0 & \mathbb{M}_2 & \dots & * \\ \vdots & \vdots & \ddots & * \\ 0 & 0 & \dots & \mathbb{M}_N \end{pmatrix},$$

and the first part of the proposition holds straightforwardly.

In order to prove (17), we fix $j < N$. Then, the points $p$ such that $h(p) = \gamma_j < 1$, lie in shifted secondary cells, let's say $\mathcal{F}_{v_1}+\delta, \dots \mathcal{F}_{v_M}+\delta$, and one can arrange the rows and columns of the matrix $\mathbb{M}_j$ such that all the points in $\mathcal{F}_{v_1}+\delta$ appear at the beginning, the points in $\mathcal{F}_{v_2}+\delta$ immediately after, and so on.

First of all, we will show that

$$\mathbb{M}_j = \begin{pmatrix} \mathbb{M}_{(j,v_1)} & 0 & \dots & 0 \\ 0 & \mathbb{M}_{(j,v_2)} & \dots & 0 \\ \vdots & \vdots & \ddots & \vdots \\ 0 & 0 & \dots & \mathbb{M}_{(j,v_M)} \end{pmatrix},$$

where $\mathbb{M}_{(j,v_k)}$ denotes the submatrix of $\mathbb{M}_j$ whose rows and columns are indexed by points in the cell $\mathcal{F}_{v_k}+\delta$.

In order to do this, consider the deformed family (15). It is straightforward to check that

$$supp\left(t^{h(p)-\omega_i(a)}x^{p-a}f_{i,\omega}\right) = (p-a, h(p)-\omega_i(a)) + Q_{i,\Omega} \subset \mathcal{Q}_\Omega + (\delta, 0)$$

where, as usual, $(i,a)$ is the row content of $p$. Moreover, the point

$$(p-a, h(p)-\omega_i(a)) = (p-a, h(p))$$

belongs to the facet of the shifted polytope $\mathcal{Q}_\Omega + (\delta,0)$ determined by and inward normal vector of the type $(v_k, l)$, $l \in \mathbb{R}_{<0}$.

Because of the last item of Lemma 3.11, there cannot be nonzero coefficients corresponding to the expansion of the polynomial

$$t^{h(p)-\omega_i(a)}x^{p-a}f_{i,\omega} \tag{19}$$

whose multidegree in $(x,t)$ lies on the boundary of $\mathcal{Q}_\Omega + (\delta,0)$ other than those in the facet determined by $(v_k, l)$. This implies that, if the point $(q,s)$ is an exponent arising in the expansion of (19) and $q$ does not belong to the shifted secondary cell $\mathcal{F}_{v_k}+\delta$, then $(q,s)$ must be an interior point of $\mathcal{Q}_\Omega + (\delta,0)$, and this implies that $s < h(q)$. Due to the remark made at the beginning of this proof, it turns out that the element indexed by $(p,q)$ in $mod(\mathbb{M}'(t))$ is zero. This gives the structure stated of $\mathbb{M}_j$.

Now, recalling Remark 3.6, and using the third item of Lemma 3.11, it turns out that $\mathbb{M}_{(j,v_k)}$ is actually a matrix of type $\mathbb{M}_{v_k,\lambda_{v_k}}$. Moreover, all matrices of this type appear in this way, so the proposition holds straightforwardly. □

Let us return to the proof of the Theorem. By the inductive hypothesis, $\det(\mathbb{M}_v^i) \neq 0$, $\forall v$, $i$, and using the special lifting $\tilde{\omega}_1, \ldots \tilde{\omega}_n, \tilde{\omega}$ in the primary cell (see Subsection 3.1.1), it is easy to see that $\det(\mathbb{M}_N) \neq 0$. Moreover,

$$\det(\mathbb{M}_N) = c_{0,\mathbf{b_0}}^{M_0} \det(\mathbb{E}_N) \neq 0, \tag{20}$$

where $\mathbb{E}_N$ is the submatrix of $\mathbb{M}_N$ formed by all the rows and columns indexed by points in non-mixed cells (see [CE1, CE2]). This proves that $\det(\mathbb{M}) \neq 0$.

Again by the inductive hypothesis, $\det(\mathbb{M}_v^i)$ equals $\mathrm{Res}_v(f_{1v}, \ldots, f_{nv})$ times a minor $\det(\mathbb{E}_v^i)$ made by choosing all rows and columns non-mixed in $\mathbb{M}_v^i$. This implies that

(21)
$$init_{\varpi} \det(\mathbb{M}) = c_{0,\mathbf{b_0}}^{M_0} \det(\mathbb{E}_N) \prod_{\mathcal{F}_v} (\mathrm{Res}(f_{1v}, \ldots, f_{nv}))^{\mathbf{d_v}} \prod_{\mathcal{F}_v} \prod_{i}^{\mathbf{d_v}} \det(\mathbb{E}_v^i).$$

By selecting all rows and columns indexed by non-mixed points in $\mathbb{M}$, a modified version of Proposition 3.12 holds for the matrix $\mathbb{E}$, which also has a block structure, and we get

$$\det(\mathbb{E}) = \det(\mathbb{E}_N) \prod_{\mathcal{F}_v} \prod_{i} \det(\mathbb{E}_v^i). \tag{22}$$

The proof of the theorem will be complete if we show that

$$init_{\varpi}(\mathrm{Res}_{\mathcal{A}}(f_0, \ldots, f_n)) = c_{0,\mathbf{b_0}}^{M_0} \prod_{\mathcal{F}_v} (\mathrm{Res}_v(f_{1v}, \ldots, f_{nv}))^{\mathbf{d_v}}.$$

Using Theorem 4.1 in [Stu2], we have that

$$init_{\varpi}(\mathrm{Res}_{\mathcal{A}}(f_0, \ldots, f_n)) = \pm \prod_{\tilde{\mathcal{F}}} (\mathrm{Res}_{\tilde{\mathcal{F}}}(f_0|_{F_0}, \ldots, f_n|_{F_n}))^{d_{\tilde{\mathcal{F}}}},$$

where $\tilde{\mathcal{F}}$ runs over all facets of the coherent mixed decomposition given by $\varpi$ in the Minkowski sum

$$\tilde{Q} := Q_0 + Q_1 + \ldots + Q_n.$$

Explicitly we have that

$$\tilde{\mathcal{F}} = F_0 + \ldots + F_n$$

with $F_i = conv(\mathcal{A}'_i),\ \mathcal{A}'_i \subset \mathcal{A}_i$,

$$f_i|_{F_i} := \sum_{a\in\mathcal{A}'_i} c_{i,a}\, x^a;$$

$\mathrm{Res}_{\tilde{\mathcal{F}}}$ is the sparse resultant associated to the data $(\mathcal{A}'_0, \ldots, \mathcal{A}'_n)$, and $d_{\tilde{\mathcal{F}}}$ equals the unique integer such that $(\mathrm{Res}_{\tilde{\mathcal{F}}}(f_0|_{F_0}, \ldots, f_n|_{F_n}))^{d_{\tilde{\mathcal{F}}}}$ has total degree $\sum_{l=0}^n \mathcal{MV}(F_0, \ldots, F_{l-1}, F_{l+1}, \ldots, F_n)$.

It is easy to see that the coherent mixed decomposition induced by $\tilde{\omega}$ over $\tilde{Q}$ is similar to the one induced over $\mathcal{Q}$. More precisely, we get a big primary cell of the form $\mathbf{b_0} + Q_1 + \ldots + Q_n$, and for each $v$ such that $\mathcal{F}_v$ is a secondary cell of $\mathcal{Q}$, there is a secondary cell $\tilde{\mathcal{F}}_v$ in $\tilde{Q}$. Moreover, the same analysis made in 3.1 says that these are all the cells of maximal dimension in $\tilde{Q}$.

For the primary cell, we get that

$$\begin{array}{rcll} f_0|_{F_0} & = & c_{0,\mathbf{b_0}}\, x^{\mathbf{b_0}}, & \\ f_i|_{F_i} & = & f_i & i = 1, \ldots, n. \end{array}$$

This implies that the unique essential set in the data $(\mathcal{A}'_0, \ldots, \mathcal{A}'_n)$ is just the singleton $\{\mathbf{b_0}\}$, so $\mathrm{Res}_{\tilde{\mathcal{F}}}(f_0|_{F_0}, \ldots, f_n|_{F_n}) = c_{0,\mathbf{b_0}}$ and $d_{\tilde{\mathcal{F}}} = \mathcal{MV}(Q_1, \ldots, Q_n) = M_0$.

Now, take a vector $v$ such that $\tilde{\mathcal{F}}_v$ is a cell of maximal dimension. This implies that $f_0|_{F_0}$ has an $n$-dimensional support and

$$f_i|_{F_i} = f_{iv} \quad i = 1, \ldots, n.$$

This implies that, in this case, the unique essential set is $\{\mathcal{A}_i \cap Q_{iv}\}_{i=1,\ldots,n}$ and $\mathrm{Res}_{\tilde{\mathcal{F}}_v}(f_0|_{F_0}, \ldots, f_n|_{F_n}) = \mathrm{Res}_v(f_{1v}, \ldots, f_{nv})$.

It remains to prove that $d_{\tilde{\mathcal{F}}_v} = \mathbf{d_v}$. To begin with, observe that

$$\mathcal{MV}(F_1, \ldots, F_n) = \mathcal{MV}(Q_{1v}, \ldots, Q_{nv}) = 0$$

due to the fact that the supports lie in a hyperplane, so the Minkowski sum of any subfamily of $\{F_1, \ldots, F_n\}$ does not have positive $n$-volume.

In order to compute the other numbers involved in the computation of $d_{\tilde{\mathcal{F}}_v}$, we shall use the recursive relation satisfied by the mixed volume ([Ber, CLO]):

$$\begin{array}{c} \mathcal{MV}_n\,(F_0, F_1, \ldots, F_{l-1}, F_{l+1}, \ldots, F_n) = \\ \sum_{v'} a_{F_0}(v')\, \mathcal{MV}'_{n-1}\,((F_1)_{v'}, (F_{l-1})_{v'}, (F_{l+1})_{v'}, \ldots, (F_n)_{v'})\,, \end{array}$$

the summation being taken over all $v'$ such that $P_{v'}$ is a facet. Here, $\mathcal{MV}'_{n-1}\,((F_i)_{v'})$ denotes the normalized mixed volume with respect to the hyperplane $v'^{\perp} \subset \mathcal{L}$ orthogonal to $v'$, and

$$a_{F_0}(v') := -\min_{m\in F_0} \langle m, v'\rangle.$$

Using the fact that $F_i = Q_{iv}$, it turns out that

$$\mathcal{M}V'_{n-1}\left((F_1)_{v'}, (F_{l-1})_{v'}, (F_{l+1})_{v'}, \dots, (F_n)_{v'}\right) = 0,$$

unless $v' = v$ or $v' = -v$. Hence, we have that

$$\begin{gathered}\mathcal{M}V_n\left(F_0, F_1, \dots, F_{l-1}, F_l, \dots, F_n\right) = \\ \left(a_{F_0}(v) + a_{F_0}(-v)\right)\mathcal{M}V'_{n-1}\left((F_1)_v, (F_{l-1})_v, (F_{l+1})_v \dots, (F_n)_v\right)\end{gathered}$$

and therefore, due to the fact that

$$\deg(\mathrm{Res}_v) = \frac{1}{\mathrm{ind}_v}\sum_{l=1}^{n}\mathcal{M}V'_{n-1}\left((F_1)_v, (F_{l-1})_v, (F_{l+1})_v \dots, (F_n)_v\right),$$

we get that

$$d_{\tilde{\mathcal{F}}_v} = \left(a_{F_0}(v) + a_{F_0}(-v)\right)\mathrm{ind}_v = \mathbf{d_v}$$

as claimed. □

**Corollary 3.13.** *For every $i = 0, 1, \dots, n$, the number of points lying in a shifted mixed cell of type $i$ is exactly $\mathcal{M}V(Q_0, \dots, Q_{i-1}, Q_{i+1}, \dots, Q_n)$.*

*Proof.* It is straightforward to check the following equalities:

$$\begin{aligned}\deg_{\mathrm{coeff}\, f_i}(\det(\mathbb{E})) &= \#\{\text{non-mixed points of type } i\}, \\ \deg_{\mathrm{coeff}\, f_i}(\det(\mathbb{M})) &= \#\{\text{non-mixed points of type } i\} \\ &\quad +\#\{\text{mixed points of type } i\}.\end{aligned}$$

□

*Remark* 3.14. We required $\mathbf{b_0}$ to be a vertex of $\mathcal{A}_0$ just in order to decrease the number of secondary cells in the algorithm, but this condition is not used in the proof. As an easy consequence, we get that

*If the family $\{f_i\}_{i=0,\dots,n}$ is essential then, every generic coefficient $c_{i,a}$ of the input support appears in $\mathrm{Res}_{\mathcal{A}}(f_0, \dots, f_n)$ with highest power $\mathcal{M}V(Q_0, \dots, Q_{i-1}, Q_{i+1}, \dots, Q_n)$.*

*Remark* 3.15. We have used the lifting algorithm of Canny and Emiris in primary cells in order to break ties, but the Theorem holds just provided that one can construct recursively Sylvester style matrices having non-zero determinant and the same degree as the resultant in the coefficients of any $f_i$. This has been already noted by Macaulay in the classical case (see [Mac, Section 6a]).

## 3.3. **Examples.**

3.3.1. *The Onedimensional Case.* Set

$$\mathcal{A}_0 := \{0, 2, 4\},\ \mathcal{A}_1 := \{4, 8\}.$$

Here, the affine lattice $\mathcal{L}$ equals $2\,\mathbb{Z}$, and the polytope $P$ is the unit segment $[0, 1]$. Set $\lambda := \frac{5}{2}$, $\delta := \frac{1}{3}$, and $\mathbf{b_0} := 0$. Then, it is straightforward to check that:

- $\mathcal{E} = [4 + \frac{1}{3}, 14 + \frac{5}{6}] \cap 2\,\mathbb{Z} = \{6, 8, 10, 12, 14\}$.
- The points lying in the shifted primary cell are $6, 8$ and $10$; the other points belong to the (unique) shifted secondary cell.
- There is a unique rule for filling the rows of the matrix $\mathbb{M}$ corresponding to the points lying in the shifted secondary cell. We have that

$$\boxed{\begin{array}{ccc} x^{12} & \mapsto & x^4\, f_1 \quad \text{(mixed)} \\ x^{14} & \mapsto & x^6\, f_1 \quad \text{(mixed).} \end{array}}$$

- Although there are infinitely many different lifting functions $\tilde{\omega}_1$, $\tilde{\omega}$ over $\mathcal{A}_1$ and $V(Q) = \{0, \frac{5}{2}\}$, they cannot produce more than two different tight coherent mixed decomposition of the segment

$$[4, 8 + \frac{5}{2}] = Q_1 + Q.$$

  Explicitly, we get the following cases:

  1.

$$\boxed{\begin{array}{cccc} x^6 & \mapsto & x^6\, f_0 & \text{(mixed)} \\ x^8 & \mapsto & x^8\, f_0 & \text{(mixed)} \\ x^{10} & \mapsto & x^2\, f_1 & \text{(non-mixed),} \end{array}}$$

  which corresponds to liftings which give the same partition than

$$\begin{array}{rcl} \tilde{\omega}_1(4, 8) & = & (0, 1) \\ \tilde{\omega}(0, \frac{5}{2}) & = & (0, 0); \end{array}$$

  2.

$$\boxed{\begin{array}{cccc} x^6 & \mapsto & x^2\, f_1 & \text{(non-mixed)} \\ x^8 & \mapsto & x^6\, f_0 & \text{(mixed)} \\ x^{10} & \mapsto & x^8\, f_0 & \text{(mixed),} \end{array}}$$

  corresponding to

$$\begin{array}{rcl} \tilde{\omega}_1(4, 8) & = & (0, 0) \\ \tilde{\omega}(0, \frac{5}{2}) & = & (0, 1). \end{array}$$

Observe that both cases give essentially the same matrix, the only difference is that the rows are indexed differently. Setting

$$\begin{array}{rcl} f_0 & = & a + bx_1^2 + cx_1^4, \\ f_1 & = & dx_1^4 + ex_1^8, \end{array} \tag{23}$$

we get the following matrix

$$\mathbb{M} := \begin{pmatrix} a & b & c & 0 & 0 \\ 0 & a & b & c & 0 \\ 0 & d & e & 0 & 0 \\ 0 & 0 & d & e & 0 \\ 0 & 0 & 0 & d & e \end{pmatrix}.$$

The matrix $\mathbb{E}$ here consists by the element $e$ which appears in the third row and column of $\mathbb{M}$. Expanding the determinant by the last column, we have that

$$\det(\mathbb{M}) = e \begin{pmatrix} a & b & c & 0 \\ 0 & a & b & c \\ 0 & d & e & 0 \\ 0 & 0 & d & e \end{pmatrix},$$

and we can easily see that the determinant of matrix on the right hand side corresponds to the Sylvester resultant for the bivariate family given by (23).

*Remark* 3.16. It is not hard to see that, for every bivariate family, every $\lambda$, $\mathbf{b_0}$ and generic $\delta$, we get a scenario similar to this example. More precisely, the algorithm given in the previous section produces the same Sylvester style matrix given by Macaulay in [Mac].

3.3.2. Let us compute the resultant of the Example 2.4. Take $\lambda = 1$, so $Q$ will be the unit square. Set also $\mathbf{b_0} = (0, 0)$. In this case, $\mathcal{L} = \mathbb{Z}^2$, the primary cell equals to 3 times $[0, 1] \times [0, 1]$ and there are two secondary cells, corresponding to the inward vectors $(-1, 0)$ and $(0, -1)$.

Setting $\delta = (\frac{2}{3}, \frac{1}{2})$, we get that $\mathcal{E}$ consists of sixteen points, nine of them lying in the shifted primary cell. Explicitly, we have that

$$\mathcal{E} = \{(a, b) \in \mathbb{Z}^2 : \ 1 \leq a, b \leq 4\}.$$

In this case, the facet polynomials $\{f_{1v}, f_{2v}\}$ can be regarded as two polynomials in one variable. So, the rows indexed by points lying in shifted secondary cells may be filled in the classical Sylvester style as in the previous example. In order to fill the rows corresponding to points lying in the shifted primary cell, we will use the following lifting functions on the ordered sets $\mathcal{A}_i = V(Q) = \{(0,0),\ (1,0)\,(0,1)\,(1,1)\}$ :

$$\begin{aligned}\tilde{\omega}_1 &= (0,1,1,2)\,;\\ \tilde{\omega}_2 &= (0,0,7,7)\,;\\ \tilde{\omega} &= (0,14,0,14)\,.\end{aligned}$$

Observe that the lifting functions are actually the restriction of a linear function on $\mathbb{R}^2$. Using any algorithm for computing convex hulls like the one given by Emiris in [Emi2], we obtain the following description for $\mathbb{M}$ :

| *row* | *coefficients of* | *cell* | *type* |
|---|---|---|---|
| $x_1x_2$ | $x_1x_2\, f_1$ | primary | non-mixed |
| $x_1^2x_2$ | $x_1^2x_2\, f_0$ | primary | mixed |
| $x_1^3x_2$ | $x_1^2x_2\, f_1$ | primary | non-mixed |
| $x_1x_2^2$ | $x_1x_2\, f_2$ | primary | non-mixed |
| $x_1^2x_2^2$ | $x_1^2x_2\, f_2$ | primary | non-mixed |
| $x_1^3x_2^2$ | $x_1^3x_2^2\, f_0$ | primary | mixed |
| $x_1x_2^3$ | $x_1x_2^2\, f_2$ | primary | non-mixed |
| $x_1^2x_2^3$ | $x_1x_2^2\, f_2$ | primary | non-mixed |
| $x_1^3x_2^3$ | $x_1^2x_2^2\, f_1$ | primary | non-mixed |
| $x_1^4x_2^3$ | $x_1^3x_2^2\, f_1$ | $(-1,0)$ − secondary | mixed |
| $x_1^4x_2^2$ | $x_1^3x_2^2\, f_2$ | $(-1,0)$ − secondary | non-mixed |
| $x_1^4x_2$ | $x_1^3x_2\, f_2$ | $(-1,0)$ − secondary | mixed |
| $x_1x_2^4$ | $x_1x_2^3\, f_2$ | $(0,-1)$ − secondary | mixed |
| $x_1^2x_2^4$ | $x_1^2x_2^3\, f_2$ | $(0,-1)$ − secondary | non-mixed |
| $x_1^3x_2^4$ | $x_1^3x_2^3\, f_2$ | $(0,-1)$ − secondary | non-mixed |
| $x_1^4x_2^4$ | $x_1^3x_2^3\, f_1$ | $(0,-1)$ − secondary | mixed |

.

Writing $f_i = a_i + b_i x_1 + c_i x_2 + d_i x_1 x_2$ and ordering the monomials as in the table, we get that

$$\mathbb{M} = \begin{pmatrix}
a_1 & b_1 & 0 & c_1 & d_1 & 0 & 0 & 0 & 0 & 0 & 0 & 0 & 0 & 0 & 0 & 0 \\
0 & a_0 & b_0 & 0 & c_0 & d_0 & 0 & 0 & 0 & 0 & 0 & 0 & 0 & 0 & 0 & 0 \\
0 & a_1 & b_1 & 0 & c_1 & d_1 & 0 & 0 & 0 & 0 & 0 & 0 & 0 & 0 & 0 & 0 \\
a_2 & b_2 & 0 & c_2 & d_2 & 0 & 0 & 0 & 0 & 0 & 0 & 0 & 0 & 0 & 0 & 0 \\
0 & a_2 & b_2 & 0 & c_2 & d_2 & 0 & 0 & 0 & 0 & 0 & 0 & 0 & 0 & 0 & 0 \\
0 & 0 & 0 & 0 & 0 & a_0 & 0 & 0 & c_0 & d_0 & b_0 & 0 & 0 & 0 & 0 & 0 \\
0 & 0 & 0 & a_2 & b_2 & 0 & c_2 & d_2 & 0 & 0 & 0 & 0 & 0 & 0 & 0 & 0 \\
0 & 0 & 0 & 0 & a_2 & b_2 & 0 & c_2 & d_2 & 0 & 0 & 0 & 0 & 0 & 0 & 0 \\
0 & 0 & 0 & 0 & a_1 & b_1 & 0 & c_1 & d_1 & 0 & 0 & 0 & 0 & 0 & 0 & 0 \\
0 & 0 & 0 & 0 & 0 & a_1 & 0 & 0 & c_1 & d_1 & b_1 & 0 & 0 & 0 & 0 & 0 \\
0 & 0 & 0 & 0 & 0 & a_2 & 0 & 0 & c_2 & d_2 & b_2 & 0 & 0 & 0 & 0 & 0 \\
0 & 0 & a_2 & 0 & 0 & c_2 & 0 & 0 & 0 & 0 & d_2 & b_2 & 0 & 0 & 0 & 0 \\
0 & 0 & 0 & 0 & 0 & 0 & a_2 & b_2 & 0 & 0 & 0 & 0 & c_2 & d_2 & 0 & 0 \\
0 & 0 & 0 & 0 & 0 & 0 & 0 & a_2 & b_2 & 0 & 0 & 0 & 0 & c_2 & d_2 & 0 \\
0 & 0 & 0 & 0 & 0 & 0 & 0 & 0 & a_2 & b_2 & 0 & 0 & 0 & 0 & c_2 & d_2 \\
0 & 0 & 0 & 0 & 0 & 0 & 0 & 0 & a_1 & b_1 & 0 & 0 & 0 & 0 & c_1 & d_1
\end{pmatrix}$$

and

$$\mathbb{E} = \begin{pmatrix}
a_1 & 0 & c_1 & d_1 & 0 & 0 & 0 & 0 & 0 & 0 \\
0 & b_1 & 0 & c_1 & 0 & 0 & 0 & 0 & 0 & 0 \\
a_2 & 0 & c_2 & d_2 & 0 & 0 & 0 & 0 & 0 & 0 \\
0 & b_2 & 0 & c_2 & 0 & 0 & 0 & 0 & 0 & 0 \\
0 & 0 & a_2 & b_2 & c_2 & d_2 & 0 & 0 & 0 & 0 \\
0 & 0 & 0 & a_2 & 0 & c_2 & d_2 & 0 & 0 & 0 \\
0 & 0 & 0 & a_1 & 0 & c_1 & d_1 & 0 & 0 & 0 \\
0 & 0 & 0 & 0 & 0 & 0 & c_2 & b_2 & 0 & 0 \\
0 & 0 & 0 & 0 & 0 & a_2 & b_2 & 0 & c_2 & d_2 \\
0 & 0 & 0 & 0 & 0 & 0 & a_2 & 0 & 0 & c_2
\end{pmatrix}$$

With the aid of *Maple*, we can check that

$$\det(\mathbb{M}) = \pm \mathrm{Res}_{\mathcal{A}}(f_0, f_1, f_n) \det(\mathbb{E}),$$

and that $\det(\mathbb{E})$ factorizes as

$$-c_2^3(-c_1 a_2 + a_1 c_2) b_2 (c_1 d_2 - d_1 c_2)(-b_2 c_1 + b_1 c_2).$$

In this easy example, one can also check that the leading term of $\det(\mathbb{M})$ as a polynomial in $a_0$ is the determinant of the following matrix:

$$
\left(\begin{array}{cccccccccccccccc}
a_1 & b_1 & 0 & c_1 & d_1 & 0 & 0 & 0 & 0 & 0 & 0 & 0 & 0 & 0 & 0 & 0 \\
0 & a_0 & 0 & 0 & 0 & 0 & 0 & 0 & 0 & 0 & 0 & 0 & 0 & 0 & 0 & 0 \\
0 & a_1 & b_1 & 0 & c_1 & d_1 & 0 & 0 & 0 & 0 & 0 & 0 & 0 & 0 & 0 & 0 \\
a_2 & b_2 & 0 & c_2 & d_2 & 0 & 0 & 0 & 0 & 0 & 0 & 0 & 0 & 0 & 0 & 0 \\
0 & a_2 & b_2 & 0 & c_2 & d_2 & 0 & 0 & 0 & 0 & 0 & 0 & 0 & 0 & 0 & 0 \\
0 & 0 & 0 & 0 & 0 & a_0 & 0 & 0 & 0 & 0 & 0 & 0 & 0 & 0 & 0 & 0 \\
0 & 0 & 0 & a_2 & b_2 & 0 & c_2 & d_2 & 0 & 0 & 0 & 0 & 0 & 0 & 0 & 0 \\
0 & 0 & 0 & 0 & a_2 & b_2 & 0 & c_2 & d_2 & 0 & 0 & 0 & 0 & 0 & 0 & 0 \\
0 & 0 & 0 & 0 & a_1 & b_1 & 0 & c_1 & d_1 & 0 & 0 & 0 & 0 & 0 & 0 & 0 \\
0 & 0 & 0 & 0 & 0 & 0 & 0 & 0 & 0 & d_1 & b_1 & 0 & 0 & 0 & 0 & 0 \\
0 & 0 & 0 & 0 & 0 & 0 & 0 & 0 & 0 & d_2 & b_2 & 0 & 0 & 0 & 0 & 0 \\
0 & 0 & 0 & 0 & 0 & 0 & 0 & 0 & 0 & 0 & d_2 & b_2 & 0 & 0 & 0 & 0 \\
0 & 0 & 0 & 0 & 0 & 0 & 0 & 0 & 0 & 0 & 0 & 0 & c_2 & d_2 & 0 & 0 \\
0 & 0 & 0 & 0 & 0 & 0 & 0 & 0 & 0 & 0 & 0 & 0 & 0 & c_2 & d_2 & 0 \\
0 & 0 & 0 & 0 & 0 & 0 & 0 & 0 & 0 & 0 & 0 & 0 & 0 & 0 & c_2 & d_2 \\
0 & 0 & 0 & 0 & 0 & 0 & 0 & 0 & 0 & 0 & 0 & 0 & 0 & 0 & c_1 & d_1
\end{array}\right) \tag{24}
$$

and we can recognize in this matrix the block structure stated in Proposition 3.12. Explicitly, we have a $9 \times 9$ big block coming from the primary cell, and two blocks arising from the secondary cells, of size $3 \times 3$ and $4 \times 4$ respectively. Computing the determinant of (24), we get:

$$-a_0^2 c_2^3(-c_1a_2 + a_1c_2)b_2(-c_2d_1 + d_2c_1)^2(b_1d_2 - b_2d_1)(c_2b_1 - c_1b_2),$$

and we can check that

$$\det(\mathbb{E}) = c_2^3(-c_1a_2 + a_1c_2)b_2(-c_2d_1 + d_2c_1)(c_2b_1 - c_1b_2),$$

and

$$
\begin{array}{rcl}
\mathrm{Res}_{(0,-1)}(c_1x_2 + d_1x_1x_2, c_2x_2 + d_2x_1x_2) & = & (-c_2d_1 + d_2c_1), \\
\mathrm{Res}_{(-1,0)}(b_1x_1 + d_1x_1x_2, b_2x_2 + d_2x_1x_2) & = & (b_1d_2 - b_2d_1)
\end{array}
$$

as expected.

3.3.3. We want to compute the sparse resultant of the family

$$
\begin{array}{rcl}
f_0 & = & a_1 + a_2x_1 + a_3x_2, \\
f_1 & = & b_1 + b_2x_1 + b_3x_2 + b_4x_1^2 + b_5x_1x_2 + b_6x_2^2, \\
f_2 & = & c_1 + c_2x_1 + c_3x_2 + c_4x_1^2 + c_5x_1x_2 + c_6x_2^2 + \\
 & & +c_7x_1^3 + c_8x_1^2x_2 + c_9x_1x_2^2 + c_{10}x_2^3.
\end{array} \tag{25}
$$

Here, $\mathcal{A}_i = \{(a,b) \in \mathbb{N}_0^2 : \, a + b \le i + 1\}, \; i = 0, 1, 2,$ the lattice $\mathcal{E}$ coincides with $\mathbb{Z}^2$ and the polytopes $Q_i$ are integer multiples of the standard simplex $S_2$. The sparse resultant coincides with the classical resultant of three homogeneous polynomials of degrees $1, 2$ and $3$ respectively, whose affinizations are the $f_i$ (see [CLO, Mac]).

In order to compute this resultant, we set $Q := 0$ and $\mathbf{b_0} := (0,0)$. The primary cell equals $5\,S_2$, and there is a unique secondary cell, corresponding to the vector $(-1,-1)$.

Setting $\delta := (\epsilon, \epsilon)$, with $1 \gg \epsilon > 0$, we have that $\mathcal{E}$ has 15 monomials, ten of them lying in the shifted primary cell. We use the following lifting function, defined on the vertices of each $Q_i$, $i = 1, 2$ ordered as follows: $\{(0,0), (k,0), (0,k)\}$ where $k = 2, 3$, and extended to the rest of the points of the input support by linearity:

$$\begin{aligned} \tilde{\omega}_1 &= (1,0,1) \\ \tilde{\omega}_2 &= (1,1,0). \end{aligned}$$

In this case, the subdivision does not depend on the value of $\tilde{\omega}$, the lifting function over the unique vertex of $Q$. As in the previous examples, the secondary cell will be filled in such a way that the facet resultants will be computed using the classical Sylvester formula for bivariate polynomials.

Explicitly, we get:

| *row* | *coefficients of* | *cell* | *type* |
|---|---|---|---|
| $x_1x_2^3$ | $x_1x_2\, f_1$ | primary | non-mixed |
| $x_1x_2^4$ | $x_1x_2^2\, f_1$ | primary | non-mixed |
| $x_1^2x_2^3$ | $x_1^2x_2\, f_1$ | primary | non-mixed |
| $x_1^4x_2$ | $x_1x_2\, f_2$ | primary | non-mixed |
| $x_1^5x_2$ | $x_1^2x_2\, f_2$ | $(-1,-1)$ – secondary | mixed |
| $x_1x_2^5$ | $x_1x_2^3\, f_1$ | $(-1,-1)$ – secondary | mixed |
| $x_1^2x_2^4$ | $x_1^2x_2^2\, f_1$ | $(-1,-1)$ – secondary | mixed |
| $x_1^3x_2^3$ | $x_1^3x_2\, f_1$ | $(-1,-1)$ – secondary | mixed |
| $x_1^4x_2^2$ | $x_1x_2^2\, f_2$ | $(-1,-1)$ – secondary | mixed |
| $x_1x_2$ | $x_1x_2\, f_0$ | primary | mixed |
| $x_1x_2^2$ | $x_1x_2^2\, f_0$ | primary | mixed |
| $x_1^2x_2$ | $x_1^2x_2\, f_0$ | primary | mixed |
| $x_1^3x_2$ | $x_1^3x_2\, f_0$ | primary | mixed |
| $x_1^3x_2^2$ | $x_1^3x_2^2\, f_0$ | primary | mixed |
| $x_1^2x_2^2$ | $x_1^2x_2^2\, f_0$ | primary | mixed |

and, indexing the matrix $\mathbb{M}$ using this order, we get

$$\begin{pmatrix}
b_6 & 0 & 0 & 0 & 0 & 0 & 0 & 0 & 0 & b_1 & b_3 & b_2 & b_5 & b_4 & 0 \\
b_3 & b_6 & b_5 & 0 & 0 & 0 & 0 & 0 & 0 & 0 & b_1 & 0 & b_2 & 0 & b_4 \\
0 & 0 & b_6 & b_4 & 0 & 0 & 0 & 0 & 0 & 0 & 0 & b_1 & b_3 & b_2 & b_5 \\
c_6 & c_{10} & c_9 & c_7 & 0 & 0 & 0 & 0 & 0 & c_1 & c_3 & c_2 & c_5 & c_4 & c_8 \\
0 & 0 & c_6 & c_4 & c_7 & 0 & c_{10} & c_9 & c_8 & 0 & 0 & c_1 & c_3 & c_2 & c_5 \\
b_1 & b_3 & b_2 & 0 & 0 & b_6 & b_5 & b_4 & 0 & 0 & 0 & 0 & 0 & 0 & 0 \\
0 & 0 & b_3 & 0 & 0 & 0 & b_6 & b_5 & b_4 & 0 & 0 & 0 & b_1 & 0 & b_2 \\
0 & 0 & 0 & b_2 & b_4 & 0 & 0 & b_6 & b_5 & 0 & 0 & 0 & 0 & b_1 & b_3 \\
c_3 & c_6 & c_5 & 0 & 0 & c_{10} & c_9 & c_8 & c_7 & 0 & c_1 & 0 & c_2 & 0 & c_4 \\
0 & 0 & 0 & 0 & 0 & 0 & 0 & 0 & 0 & a_1 & a_3 & a_2 & 0 & 0 & 0 \\
a_3 & 0 & 0 & 0 & 0 & 0 & 0 & 0 & 0 & 0 & a_1 & 0 & a_2 & 0 & 0 \\
0 & 0 & 0 & 0 & 0 & 0 & 0 & 0 & 0 & 0 & 0 & a_1 & a_3 & a_2 & 0 \\
0 & 0 & 0 & a_2 & 0 & 0 & 0 & 0 & 0 & 0 & 0 & 0 & 0 & a_1 & a_3 \\
0 & 0 & 0 & 0 & 0 & 0 & 0 & a_3 & a_2 & 0 & 0 & 0 & 0 & 0 & a_1 \\
0 & 0 & a_3 & 0 & 0 & 0 & 0 & 0 & 0 & 0 & 0 & 0 & a_1 & 0 & a_2
\end{pmatrix}$$

In this case, the extraneous factor $\mathbb{E}$ is the submatrix made by choosing the first four rows and columns of $\mathbb{M}$, and its determinant equals

$$b_6\,(c_7\,b_6^2 - b_4\,b_6\,c_9 + b_4\,b_5\,c_{10}).$$

Observe that all points lying in shifted non-mixed cells are actually in the shifted the primary cell, so there is a priori no significant block structure in this matrix.

*Curiosity:* The matrix constructed here is exactly the same given by Canny and Emiris in the last section of [CE2] (see also section 3.1.4 of [Emi1]) in order to show that, in their construction, the extraneous factor is not always the determinant of the minor formed by choosing all rows and columns indexed by non-mixed points. Of course, they work with another definition of non-mixed cells!

To be more precise, their algorithm is not recursive. They produce this matrix by applying a lifting algorithm to the polytopes $Q_0$, $Q_1$ and $Q_2$ using the linear functions $l_0 := 10^4x_1 + 10^3x_2,\ l_1 := 10^5x_1,\ l_2 := 10^2x_1+10x_2$. By taking the lower hull of the lifted polytopes $(Q_i, l_i(Q_i))$, a tight coherent mixed decomposition of the Newton polytope $Q_0+Q_1+Q_2$ holds. Using $\delta$ as before, the same matrix is constructed but the points lying in the shifted non-mixed cells are in correspondence with the monomials

$$\{x_1x_2^3, x_1x_2^4, x_1^2x_2^3, x_1^5x_2\}.$$

By taking the determinant of the submatrix of $\mathbb{M}$ made choosing the rows and columns indexed by these monomials, we get $b_6^3\,c_7$.

3.4. **An Overview of Macaulay's Classical Formulas.** In this section, we will see how the formulas given by Macaulay in [Mac] can be recovered with our methods. In order to have a notation similar to Macaulay's original paper, we shall deal with $n$ generic polynomials in $n-1$ variables $x_1, \dots, x_{n-1}$ of total degree less or equal than $m_1, \dots, m_n$ respectively. In our terminology, the input supports are integer multiples of the standard simplex $S_{n-1}$. More precisely, the polynomials will be denoted as $C_1$, $C_2, \dots$ , $C_n$ and

$$supp(C_i) = m_i\, S_{n-1}.$$

Actually, Macaulay worked with $n$ homogeneous polynomials in $n$ variables, but it turns out that the homogeneous resultant of these polynomials coincides with the sparse resultant associated to the supports $C_i$ as in Section 2.

*Remark* 3.17. Setting $t_n := \sum_{i=1}^{n}(m_i - 1)$, and using the algorithm given by Canny and Emiris, it is also possible to recover Macaulay's classical formula in degree $t = t_n + 1$. Moreover, the extraneous factor in Macaulay's original formulation is exactly the minor formed by using all rows and columns indexed by points lying in shifted non-mixed cells which they get with their methods (see [CE2] and [Emi1]).

In section 3 of [Mac], Macaulay constructed a Sylvester style matrix whose determinant is denoted by $D(n, t)$, for every $t \in \mathbb{N}_{\geq 0}$. The matrix has its rows and columns indexed by all monomials of total degree less or equal than $t$. For every $t > t_n$, it turns out that $D(n, t)$ is a nonzero multiple of $\mathrm{Res}_{\mathcal{A}}(C_1, \dots, C_n)$ denoted in that paper as $R(n, t)$. [2]

His construction is recursive in the following sense: let $\nu$ be the inward normal vector $(-1, -1, \dots, -1)$; it is easy to see that

$$C_{1\nu},\, C_{2\nu}, \dots, C_{n-1\nu}. \tag{26}$$

is a family of $n-1$ homogeneous polynomials in the variables $x_1, \dots, x_{n-1}$, so after setting $x_{n-1}$ equal to one, it may be regarded as a sparse family with support in

$$m_1\, S_{n-2}, m_2\, S_{n-2}, \dots, m_{n-1}\, S_{n-2},$$

so $D(n-1, j)$ will be the determinant of a Sylvester matrix made with the same rules as $D(n, j)$ but with $n-2$ variables, using the polynomials (26).

[2]For lower values of $t$, Macaulay also proposed a matrix with similar properties, but not of Sylvester type.

In section 5 of [Mac], he established the following theorem:

$$\frac{D(n,t)}{R(n,t)} = \prod_{j=0}^{m_n-1} \frac{D(n-1,t-j)}{R(n-1,t-j)} \prod_{k=1}^{t-m_n} D(n-1,k), \tag{27}$$

which is the cornerstone of the main result given in section 6, namely, that $R(n,t)$ can be recovered as the quotient of $D(n,t)$ by the minor of obtained by omitting all rows and columns corresponding to monomials *reduced* in all the variables.

Here, a monomial $x^\alpha = x_1^{\alpha_1} \dots, x_{n-1}^{\alpha_{n-1}}$ is said to be reduced if there exists a unique $i \in \{1, \dots, n\}$ such that $x_i^{m_i}$ divides $x^\alpha\, x_n^{j-|\alpha|}$.

Macaulay gave also a recursive structure of the extraneous factor (in his notation, $\Delta(n,t)$) as follows:

$$\Delta(n,t) = \prod_{j=0}^{m_n-1} \Delta(n-1,t-j) \prod_{k=1}^{t-m_n} D(n-1,k) \tag{28}$$

(see section 6 of [Mac]).

In order to see this construction in light of the results presented in the previous section, extra care must be taken, because Macaulay's construction produces a determinant which has the same degree as the resultant in the coefficients of $C_n$, so we shall modify our algorithm in such a way that the role of $f_0$ is played by $C_n$; similarly, the role of $f_1$ will be played by $C_{n-1}$, and so on.

The polytope $P$ will be the standard simplex $S_{n-1}$. Given $t > t_n$, in order to define $Q$, set

$$\lambda := t - t_n - 1,$$

so $\mathcal{Q}$ will be equal to

$$(m_1 + \ldots + m_n + \lambda) S_{n-1} = \{(m_1, \dots, m_{n-1}) \in \mathbb{R}^{n-1} : \ 0 \leq m_i \leq t+n-1\},$$

and $\mathcal{L} = \mathbb{Z}^{n-1}$. By taking $\delta = (\epsilon, \dots, \epsilon)$, with $1 \gg \epsilon > 0$, we get that

$$\mathcal{E} = \{(\alpha_1, \dots, \alpha_{n-1}) \in \mathbb{Z}^{n-1} : \ 1 \leq \alpha_i, \ \sum_i \alpha_i \leq t+n-1\}.$$

Setting

$$\begin{array}{rcl} \mathcal{E} & \to & \{\beta \in \mathbb{N}_0^{n-1} : \sum \beta_i \leq t\} \\ (\alpha_1, \dots, \alpha_{n-1}) & \mapsto & (\alpha_1 - 1, \dots, \alpha_{n-1} - 1), \end{array} \tag{29}$$

we get a bijection between our support and the one used by Macaulay for computing the matrix whose determinant is $D(n,t)$. By choosing $\mathbf{b_0} = (0, \dots, 0)$, we can check that there is only one secondary cell, associated to the vector $\nu$, whose $\nu-$diameter is exactly $m_n$.

Using lifting functions $\tilde{\omega_1}, \dots, \tilde{\omega_{n-1}}, \tilde{\omega}$ as in section 8 of [CE2], it is possible to get a subdivision of the primary cell such that the points lying in the 0−mixed cells are bijectively associated (with the bijection given in (29)) with those reduced monomials which are divisible by $x_n^{m_n}$. Lifting functions with the same properties should be used in the recursive steps.

Now, comparing equations (28) and (22), it is not hard to check that

$$\prod_{j=0}^{m_n-1} \Delta(n-1, t-j) = \prod_{\mathcal{F}_v} \prod_i \det(\mathbb{E}_v^i).$$

Hence,

$$\det(\mathbb{E}_N) = \prod_{k=1}^{t-m_n} D(n-1, k).$$

## 4. The general case

We will extend here the results of the previous section by assuming only that the family $\{\mathcal{A}_i\}_{0\leq i\leq n}$ is essential, without any other condition on the supports. Observe that this hypothesis ensures that $\mathrm{Res}_{\mathcal{A}}(f_0, \ldots, f_n)$ is non trivial.

As before, set

$$Q_i = conv(\mathcal{A}_i), \quad i = 0, \ldots, n,$$

and consider the Minkowski sum $Q_0 + \ldots + Q_n$. As the family of supports is essential, this polytope must be $n$-dimensional. So, $\mathcal{L}_{\mathbb{R}} = \mathbb{R}^n$.

In order to give more generality to our algorithm, let us consider a "generic" polytope $Q$ (in a sense which will be determined later, see for instance Remark 4.4), with vertices in $\mathcal{L} \otimes \mathbb{Q}$. The algorithm will produce a Sylvester style matrix whose rows and columns will be indexed by the elements of

$$\mathcal{E} := (Q_0 + Q_1 + \ldots + Q_n + Q) \cap \mathcal{L}, \tag{30}$$

and whose determinant will be a non zero multiple of the sparse resultant. The extraneous factor will be a minor of this matrix.

As before, the algorithm will be recursive on the dimension of the polytope

$$\mathcal{Q} := Q_0 + Q_1 + \ldots + Q_n + Q,$$

and in its intermediate steps will use the subdivision technique of Canny and Emiris (see [CE1, CE2, Stu2]).

4.1. **Construction of the Sylvester matrix.** Let $Q$ be polytope as before, $V(Q) \subset \mathcal{L} \otimes \mathbb{Q}$ be the set of vertices of $Q$.

Let us pick a vertex $\mathbf{b_0} \in \mathcal{A}_0$, and consider the same lifting functions $\omega_i$, $\omega$, defined in (4). Consider also

$$\begin{array}{rcl} Q_{i,\Omega} & := & conv\{(a, \omega_i(a)) : a \in \mathcal{A}_i\}, \\ Q_\Omega & := & conv\{(b, \omega(b)) : b \in V(Q)\}, \end{array}$$

and the coherent mixed decomposition $\Delta_{i,\Omega}$ (resp. $\Delta_\Omega$) of the polytopes $Q_i$ (resp. $Q$) given by projecting the upper envelope of $Q_{i,\Omega}$ (resp. $Q_\Omega$).

In the decomposition of $\mathcal{Q}$, each cell is of the form

$$F = F_0 + F_1 + \ldots + F_n + F,$$

where $F_i$ (resp. $F$) is a cell in $\Delta_{i,\Omega}$, (resp. $\Delta_\Omega$).

Our lifting functions are again mostly trivial. So, we can compute the cells which will appear in the subdivision: in $\Delta_{1,\Omega}, \ldots, \Delta_{n,\Omega}, \Delta_\Omega$ we may have two types of cells as in the previous section.

The decomposition $\Delta_{0,\Omega}$ has the following cells: the face determined by $(0, -1)$ is just the singleton $\{\mathbf{b_0}\}$, and it is straightforward to check that, for every $k \geq 0$, every $k$-dimensional cell in the decomposition is a $k$-dimensional face of $Q_0$ which contains $\{\mathbf{b_0}\}$, or the convex hull of $\mathbf{b_0}$ and a $(k-1)$-dimensional face of $Q_0$ which does not contain this point. We shall denote it by $F_{0,v}$, where $v$ is some interior primitive normal vector of that face.

So, all cells in the decomposition of $\mathcal{Q}$ are as follows:

- the *primary cell:* $\{\mathbf{b_0}\} + Q_1 + \ldots + Q_n + Q$. Observe that it has always dimension $n$ due to the fact that the family $\{\mathcal{A}_i\}_{i=0,1,\cdots,n}$ is essential.
- $F_{0,v} + Q_{1v} + \ldots + Q_{nv} + Q_v$ for some $v \in \mathbb{R}^n$. Those being $n$-dimensional will be called *secondary*, and they will be uniquely determined by their primitive inward vector $v \in \mathbb{R}^n$.

As in the previous section, the Sylvester matrix will take into account if the points lie in the primary cell or not. The first assumption we will make on $Q$ is that all points in $\mathcal{Q} \cap \mathcal{L}$ must belong to the interior of a maximal cell (primary or secondary).

4.1.1. *Points in the primary cell.* Proceed exactly as in 3.1.1: choose generic lifting functions $\tilde{\omega}_1, \ldots, \tilde{\omega}_n, \tilde{\omega}$ over en $\mathcal{A}_1, \ldots, \mathcal{A}_n$, $V(Q)$ respectively, in such a way that they produce a tight coherent mixed decomposition of the Minkowski sum

$$Q_1 + \ldots + Q_n + Q$$

by taking the upper envelope. Every $n$-dimensional cell $\mathcal{F}$ in the decomposition may be obtained as

$$F_1 + F_2 + \ldots + F_n + F,$$

where $F_i$ is a cell of $\Delta_{i,\tilde{\omega}}$, $F$ is a cell of $\Delta_{\tilde{\omega}}$,

$$n = \dim(F_1) + \ldots + \dim(F_n) + \dim(F)$$

and at least one of these dimensions is zero.

The concepts of row content, mixed cells of type 0, and the entries of $\mathbb{M}$ whose row coordinates are indexed by

$$p \in (\{\mathbf{b_0}\} + Q_1 + \ldots + Q_n + Q) \cap \mathcal{L},$$

are defined in the same way as in 3.1.1.

4.1.2. *Points in the secondary cell.* Here, we cannot use the recursive step given in 3.1.2, so we must proceed with some care. Let $v \in \mathbb{R}^n$ the primitive inward normal vector of the $n$-dimensional cell

$$\mathcal{F}_v = F_{0,v} + Q_{1_v} + \ldots + Q_{n_v} + Q_v. \tag{31}$$

Consider the family of polynomials associated to the $v$-facet

$$f_{i_v}(x_1, \ldots, x_n) = \sum_{a \in Q_{i_v} \cap \mathcal{A}_i} c_{i,a}\, x^a \ (i = 1, \ldots, n). \tag{32}$$

We cannot claim now that the family

$$\{Q_{i_v} \cap \mathcal{A}_i\}_{1 \leq i \leq n} \tag{33}$$

is essential, but certainly there always exists a subfamily of indices $1 \leq i_1 < i_2 < \cdots < i_k \leq n$, such that

$$\{Q_{i_{j\,v}} \cap \mathcal{A}_{i_j}\} \tag{34}$$

is essential. It is not always true that (34) is the unique essential subfamily of (33), but we can assume w.l.o.g. that there exists an essential subfamily which is minimal with respect to the inclusion. Moreover, we may suppose that this subfamily is indexed as follows:

$$\{Q_{i_v} \cap \mathcal{A}_i\}_{1 \leq i \leq k}. \tag{35}$$

This means that the sparse resultant of the polynomials $f_{1_v}, f_{2_v}, \cdots, f_{k_v}$ with respect to the supports (35) is non trivial. We shall denote it as

$$\mathrm{Res}_v(f_{1_v}, \ldots, f_{k_v}). \tag{36}$$

In order to mimic the inductive step of the previous section, suppose w.l.o.g. that $0 \in \mathcal{A}_{iv}$, $i = 1, \cdots, n$. Then, we can consider the chain of lattices

$$\mathcal{L}_{\mathcal{A}_{1_v} + \cdots + \mathcal{A}_{k_v}} \subset \mathcal{L}_{\mathcal{A}_{1_v} + \cdots + \mathcal{A}_{n_v}} \subset \mathcal{L}_v,$$

and we know that $\mathcal{L}_{\mathcal{A}_{1_v}+\cdots+\mathcal{A}_{k_v}}$ has dimension equal to $k-1$. Hence, we have that

$$Q_{1_v} + Q_{2_v} + \cdots + Q_{k_v} \subset \mathcal{L}_{\mathcal{A}_{1_v}+\cdots+\mathcal{A}_{k_v}} \otimes \mathbb{R}. \tag{37}$$

Set $\tilde{Q}_v := F_{0,v} + Q_{k+1_v} + \ldots + Q_{n_v} + Q_v$. We have that

$$(F_{0,v} + Q_{1_v} + \ldots + Q_{n_v} + Q_v) \cap \mathcal{L} = \left(Q_{1_v} + \ldots + Q_{k_v} + \tilde{Q}_v\right) \cap \mathcal{L}. \tag{38}$$

For a sublattice $G \subset \mathcal{L}$, its saturation will be denoted by $s(G)$. Consider the orthogonal decomposition $\mathcal{L} = s(\mathcal{L}_{\mathcal{A}_{1_v}+\cdots+\mathcal{A}_{k_v}}) \oplus \mathcal{L}'$. Every point $p \in \mathcal{L}$ may be written as $p = \mathrm{q}_j + p_{kv} + p^v$, where $\mathrm{q}_1, \ldots, \mathrm{q}_{\mathrm{ind}_v}$ are coset representatives for $\mathcal{L}_{\mathcal{A}_{1_v}+\cdots+\mathcal{A}_{k_v}}$ in $s(\mathcal{L}_{\mathcal{A}_{1_v}+\cdots+\mathcal{A}_{k_v}})$, the number of such cosets is $\mathrm{ind}_v$, $p_{kv} \in \mathcal{L}_{\mathcal{A}_{1_v}+\cdots+\mathcal{A}_{k_v}}$ and $p^v \in \mathcal{L}'$.

As (38) is finite, we can intersect $\tilde{Q}_v$ with finitely many rational translates of $\mathcal{L}_{\mathcal{A}_{1_v}+\cdots+\mathcal{A}_{k_v}} \otimes \mathbb{R}$, and get polytopes

$$\tilde{Q}_v^m := (\tilde{q}_{mv} + \mathcal{L}_{\mathcal{A}_{1v}+\ldots+\mathcal{A}_{kv}} \otimes \mathbb{R}) \cap \tilde{Q}_v, \quad m = 1, ..., M,$$

such that (38) is equal to the disjoint union of the sets

$$\left(Q_{1_v} + \ldots + Q_{k_v} + \tilde{Q}_v^m\right) \cap \left(\tilde{q}_{mv} + \mathrm{q}_j + \mathcal{L}_{\mathcal{A}_{1_v}+\cdots+\mathcal{A}_{k_v}}\right),$$

where $m = 1, \ldots, M$, and$j = 1, \ldots, \mathrm{ind}_v$. Now, it is easy to see that $\tilde{\mathbf{Q}}_v^m := -\mathrm{q}_j - \tilde{q}_{mv} + \tilde{Q}_v^m$ is a polytope with vertices in $\mathcal{L}_{\mathcal{A}_{1_v}+\cdots+\mathcal{A}_{k_v}} \otimes \mathbb{Q}$. So, as in the previous section, for every $p$ in (38), there exists a unique $m \in \{1, \ldots, M\}$ such that $\mathbf{p} := p - \tilde{q}_{mv} - \mathrm{q}_j$ belongs to

$$\left(Q_{1_v} + \ldots + Q_{k_v} + \tilde{\mathbf{Q}}_v^m\right) \cap \mathcal{L}_{\mathcal{A}_{1_v}+\cdots+\mathcal{A}_{k_v}}. \tag{39}$$

Consider also the family

$$\{f_{i_v}(x_1, \ldots, x_n)\}_{1 \leq i \leq k}. \tag{40}$$

We may now apply the inductive hypothesis to (39) and construct, as in the previous section square matrices $\mathbb{M}_{v,\tilde{\mathbf{Q}}_v^m}$ indexed by the points in (38), whose determinants are non zero multiples of $\mathrm{Res}_v(f_{1_v}, \ldots, f_{kv})$, each of them with the same degree in the coefficients of $f_{1_v}$ as the resultant.

As in the primary cell, each of these matrices will have, in the row indexed by the point $p$, the coordinates, in the monomial basis, of the expansion of $x^{p-a} f_{iv}(x)$ for some pair $(a, i)$. In order to define the entries of the matrix $\mathbb{M}$, we shall proceed as before: the entry indexed by $(p, p')$ wil be the coefficient of $x^{p'}$ in the expansion of $x^{p-a} f_i(x)$.

Now that we have the matrix well defined, it remains to decide which points lying in secondary cells will be mixed and which not. In order to do this, we shall begin with the following

**Definition 4.1.** We shall say that the vector $v$ is *admissible* if there exists a unique essential subfamily of (33).

Recall that Theorem 4.1 in [Stu2] applied to the weight $\overline{\omega}$ defined in (13) gives

$$init_{\overline{\omega}}(\mathrm{Res}_{\mathcal{A}}(f_0,\dots,f_n)) = \pm \prod_{\tilde{\mathcal{F}}} \left(\mathrm{Res}_{\tilde{\mathcal{F}}}(f_0|_{F_0},\dots,f_n|_{F_n})\right)^{d_{\tilde{\mathcal{F}}}},$$

where $\tilde{\mathcal{F}}$ ranges over all facets of the coherent mixed decomposition given by $\overline{\omega}$ on $\tilde{Q} := Q_0 + Q_1 + \dots + Q_n$,

$$\begin{aligned}
&\tilde{\mathcal{F}} = F_0 + \dots + F_n,\\
&F_i = conv(\mathcal{A}'_i),\ \mathcal{A}'_i \subset \mathcal{A}_i,\\
&f_i|_{F_i} := \textstyle\sum_{a\in\mathcal{A}'_i} c_{i,a}\, x^a,
\end{aligned}$$

$\mathrm{Res}_{\tilde{\mathcal{F}}}$ being the sparse resultant associated with the family of supports $(\mathcal{A}'_0,\dots,\mathcal{A}'_n)$, and $d_{\tilde{\mathcal{F}}}$ equals the unique integer such that

$$\left(\mathrm{Res}_{\tilde{\mathcal{F}}}(f_0|_{F_0},\dots,f_n|_{F_n})\right)^{d_{\tilde{\mathcal{F}}}}$$

has total degree $\sum_{l=0}^{n} \mathcal{M}V(F_0,\dots,F_{l-1},F_{l+1},\dots,F_n)$.

*Remark* 4.2. Observe that $\mathrm{Res}_{\tilde{\mathcal{F}}} \neq 1$ if and only if $\tilde{\mathcal{F}}$ is associated with an admissible vector $v$.

The integer $d_{\tilde{\mathcal{F}}}$ will allow us to "choose" the mixed points as follows:

- if $v$ is not admissible, then *all* points lying in $\mathcal{F}_v \cap \mathcal{L}$ will be non mixed.
- If $v$ is admissible, then we may choose $d_{\tilde{\mathcal{F}}_v}$ of the matrices $\mathbb{M}_{v,\tilde{Q}}$ made by subdividing $\mathcal{F}_v \cap \mathcal{L}$ in sets of the form (39). We will see in Proposition 4.3 that there is at least this number of such matrices.
  - If a point indexes one of the $d_{\tilde{\mathcal{F}}_v}$ matrices $\mathbb{M}_{v,\tilde{Q}}$ and is "mixed" for that matrix, then it will be mixed for the matrix $\mathbb{M}$.
  - All the other points lying in $\mathcal{F}_v \cap \mathcal{L}$ will be non mixed.

Again, it is easy to see that, in the case $n = 1$, we get the classical Sylvester style matrices([Syl, Mac]).

**Proposition 4.3.** *If $v$ is admissible, then the number of sets of the form (39) is greater or equal than $d_{\tilde{\mathcal{F}}_v}$.*

*Proof of the Proposition:* It is easy to see that $init_{\overline{\omega}}(\mathrm{Res}_{\mathcal{A}}(f_0,\dots,f_n))$ is actually the sparse resultant $\mathrm{Res}_{\mathcal{A}}$ specialized in the family

$$c_{0,\mathbf{b_0}}x^{\mathbf{b_0}}, f_1, \cdots, f_n.$$

In order to have a nice interpretation of $d_{\tilde{\mathcal{F}}_{v_1}}$, we will use Minimair's formula given in [Min, Theorem 1] for computing $\mathrm{Res}_{\mathcal{A}}(c_{0,\mathbf{b_0}}x^{\mathbf{b_0}}, f_1, \cdots, f_n)$. First of all, observe that the hypothesis of Minimair's Theorem is satisfied, due to the fact that the family $\mathcal{A}_0, \cdots, \mathcal{A}_n$ is essential. This implies that the unique essential subfamily of $\{\mathbf{b_0}\}, \mathcal{A}_1, \cdots, \mathcal{A}_n$ is $\{\mathbf{b_0}\}$.

Using this formula, we have that:

$$\mathrm{Res}_{\mathcal{A}}(c_{0,\mathbf{b_0}}x^{\mathbf{b_0}}, f_1, \cdots, f_n) = c_{0,\mathbf{b_0}}^{M_0}\prod_v \mathrm{Res}_v(f_{i_{1v}}, \cdots, f_{i_{kv}})^{e_v},$$

where the product ranges over all primitive inward normal vectors of the facets of the Newton polytope of $\mathcal{A}_1 + \cdots + \mathcal{A}_n$, and

$$e_v = (a_{F_0}(v) + a_{F_0}(-v))\ [\mathcal{L}_v : \mathcal{L}_{\mathcal{A}_{1v}+\cdots+\mathcal{A}_{nv}}]\ \mathbf{e}'_{\mathbf{v}}, \tag{41}$$

- $F_0$ is the convex hull of $\{\mathbf{b_0}\} \cup \mathcal{A}_{0v}$,
- If $\{A_{i_{1v}}, \cdots A_{i_{kv}}\}$ is the unique essential subfamily of $\{\mathcal{A}_{1v}, \cdots, \mathcal{A}_{nv}\}$, then $\mathbf{e}'_{\mathbf{v}}$ is defined as follows:
  1. Let $\mathcal{L}_{\mathcal{A}_{1v}+\cdots+\mathcal{A}_{nv}} = s(\mathcal{L}_{A_{i_{1v}}+\cdots+A_{i_{kv}}}) \oplus \mathcal{L}^{\circ}$ be the orthogonal decomposition, and denote by $\pi$ the projection over the second factor.
  2. Define $\mathbf{e}'_{\mathbf{v}} := \mathcal{MV}\left(\pi(Q_{iv})\right)_{i\notin\{i_1,\cdots,i_k\}}$, where $\mathcal{MV}(\cdot)$ denotes the normalized mixed volume with respect to $\mathcal{L}^{\circ}$.

Now we will prove the proposition. We can suppose w.l.o.g. that $Q = 0$, because adding a polytope to the Minkowski sum cannot make decrease the number of sets of the form (39). Moreover, we may suppose w.l.o.g. that the unique essential family is $\mathcal{A}_{1v}, \cdots, \mathcal{A}_{kv}$.

If $v$ is an admissible vector, then $e_v$ and $d_{\mathcal{F}_v}$ coincide. Following the construction of the Sylvester style matrix for the generalized unmixed case (Section 3.1), it turns out that $\mathcal{F}_v \cap \mathcal{L}$ may be decomposed as a disjoint union of $(a_{F_0}(v) + a_{F_0}(-v))\,[\mathcal{L}_v : \mathcal{L}_{\mathcal{A}_{1v}+\cdots+\mathcal{A}_{nv}}]$ sets parallel to $Q_{1v} + \cdots + Q_{nv}$.

From each of these partitions we may have at least as many subsets of the form (39) as integer points are in

$$\pi(Q_{k+1_v}) + \cdots + \pi(Q_{nv}). \tag{42}$$

It is well-known that the mixed volume of these polytopes may be computed as the number of integer points lying in a subset of a displacement of (42) (see for instance [Emi1, HS]), so $\mathbf{e}'_{\mathbf{v}}$ will be less or equal than the number of integer points in (42). $\square$

*Remark* 4.4. The reader can check that the role played by $\delta$ in the previous section is played here by the additional polytope $Q$ : at each step of the recursion, we need to impose some conditions on the different additional polytopes of lower dimensions in order to guarantee that all integer points are in the interior of a cell. It is easy to see that this happens for a "generic" polytope. One can set, for instance,

$$Q = \text{fixed polytope} + \delta,$$

with $\delta$ generic as in the previous section.

4.2. **Generalized Macaulay style formula.** What follows may be regarded as the main result of this paper, and an extension of Theorem 3.8:

**Theorem 4.5.** $\mathbb{M}$ *is a Sylvester style matrix and* $\det(\mathbb{M}) \neq 0$. *Moreover, we have the following formula à la Macaulay:*

$$\det(\mathbb{M}) = \mathrm{Res}_{\mathcal{A}}(f_0, \ldots, f_n) \det(\mathbb{E}),$$

*where* $\mathbb{E}$ *is the square submatrix of* $\mathbb{M}$ *made by omitting all rows and columns indexed by non mixed points.*

As $\mathbb{E}$ does not contain coefficients of $f_0$, we get again that $\det(\mathbb{M})$ has the same degree as $\mathrm{Res}_{\mathcal{A}}(f_0, \ldots, f_n)$ in the coefficients of $f_0$. Replacing the rol played by $f_0$ with $f_i$, $i = 1, \cdots, n$, we have a formula for computing $\mathrm{Res}_{\mathcal{A}}(f_0, \ldots, f_n)$ as the gcd of $n+1$ determinants.

*Proof.* The same argument given in the proof of Theorem 3.8 may be applied to this situation in order to see that $\mathbb{M}$ is a Sylvester style matrix. So, we have that $\mathrm{Res}_{\mathcal{A}}(f_0, \ldots, f_n)$ divides to $\det(\mathbb{M})$.

As before, we shall prove that $\det(\mathbb{M})$ is not identically zero by showing that its initial term with respect to $\overline{\omega}$ is non zero. Again we have

$$\deg_{\text{coeff }(f_0)}(\mathrm{Res}_{\mathcal{A}}(f_0, \ldots, f_n)) = \deg_{\text{coeff }(f_0)}(\det(\mathbb{M})).$$

Then,

$$\frac{\det(\mathbb{M})}{\mathrm{Res}_{\mathcal{A}}(f_0, \ldots, f_n)} = \frac{init_{\overline{\omega}}(\det(\mathbb{M}))}{init_{\overline{\omega}}(\mathrm{Res}_{\mathcal{A}}(f_0, \ldots, f_n))}$$
$$= \frac{\text{coeffficient of } c_{0,\mathbf{b_0}}^{M_0} \text{ in } \det(\mathbb{M})}{\text{coefficient of } c_{0,\mathbf{b_0}}^{M_0} \text{ in } \mathrm{Res}_{\mathcal{A}}(f_0, \ldots, f_n)}.$$

We will see that the numerator of this fraction is non zero, and that the ratio is $\det(\mathbb{E})$ by induction on $n$. The initial case was already covered in the previous section. Suppose $n > 1$, and let us introduce again a parameter of deformation $t$ :

$$\begin{array}{rcl} f_{0,\omega} & = & c_{0,\mathbf{b_0}}\, t\, x^{\mathbf{b_0}} + \sum_{a\in\mathcal{A}_i\setminus\{\mathbf{b_0}\}} c_{i,a}x^a, \\ f_{i,\omega} & = & f_i \qquad\qquad\qquad\qquad i \geq 1. \end{array}$$

Consider the modified matrix $\mathbb{M}\left(c_{i,a}\, t^{\omega_i(a)}\right)$. For $p \in \mathcal{E}$, let $h(p)$ be the largest rational number such that

$$(p, h(p)) \in \mathcal{Q}_\Omega = Q_{0,\Omega} + Q_{1,\Omega} + \ldots + Q_{n,\Omega} + Q_\Omega.$$

For every $p \in \mathcal{E}$, we shall multiply all the entries in the row indexed by $p$ by $t^{h(p)-\omega_i(a)}$ where -as usual- $(i, a)$ denotes the row content of $p$. We shall denote this matrix with $\mathbb{M}'(t)$. It is not hard to see that the leading coefficient of $\det(\mathbb{M}'(t))$ is $init_{\overline{\omega}}(\det(\mathbb{M}))$.

The following assertions may be proven mutatis mutandis the results given in the previous section:

**Lemma 4.6.** *It turns out that*

1. $0 < h(p) \leq 1$, *for every* $p \in \mathcal{E}$.
2. $h(p) = 1$ *if and only if* $p$ *belongs to the primary cell.*
3. *If* $p$ *and* $q$ *both belong to the same secondary cell, let us say* $\mathcal{F}_v$, *then*
$$h(p) = h(q) \iff \langle p, v\rangle = \langle q, v\rangle.$$
4. *If* $p \in \mathcal{F}_v \cap \mathcal{L}$ *has row content* $(i, a)$ *and* $v' \neq \mu\, v,\ \mu > 0$, *then*

$$((p - a, h(p)) + Q_{i,\Omega}) \cap \mathcal{Q}_{\Omega(v',l')} = \emptyset,\ \forall\, l' \in \mathbb{R}_{<0}. \tag{43}$$

**Proposition 4.7.** *Let* $0 < \gamma_1 < \gamma_2 < \ldots < \gamma_N = 1$ *be the different values of* $h(p)$ *for* $p \in \mathbb{E}$. *Then, the leading coefficient of* $\det(\mathbb{M}'(t))$ *(regarded as a polynomial in* $t$*) factorizes as follows:*

$$\prod_{j=1}^{N} \det\left(\mathbb{M}_j\right), \tag{44}$$

*where* $\mathbb{M}_j$ *is the square submatrix of* $\mathbb{M}$ *made by choosing all rows and columns indexed by those points* $p$ *such that* $h(p) = \gamma_j$.

*This product may be also factorized as follows:*

$$\det(\mathbb{M}_N) \prod_{\mathcal{F}_v} \left(\det(\mathbb{M}_v^1) \ldots \det(\mathbb{M}_v^{\mathbf{d_v}})\right), \tag{45}$$

*donde*

- *the product ranges over all secondary cells* $\mathcal{F}_v$,
- *For every* $v$, $\mathbf{d_v}$ *is the number of sets of the form (39) obtained by subdividing* $\mathcal{F}_v$,

- $\mathbb{M}_v^i$ *is a matrix of the type* $\mathbb{M}_{v,\tilde{Q}}$, *the latter being defined just after* (40).

By the inductive hypothesis, it turns out that $\det(\mathbb{M}_v^i) \neq 0,\ \forall v,\ i$. Using the same argument given in the proof of Theorem 3.8, we also get that

$$\det(\mathbb{M}_N) = c_{0,\mathbf{b_0}}^{M_0} \det(\mathbb{E}_N) \neq 0.$$

Here, $\mathbb{E}_N$ is the square submatrix of $\mathbb{M}_N$ made by all those rows and columns indexed by non mixed points. This implies that $\det(\mathbb{M}) \neq 0$. Moreover,

(46)

$$init_{\varpi} \det(\mathbb{M}) = c_{0,\mathbf{b_0}}^{M_0} \det(\mathbb{E}_N) \prod_{\mathcal{F}_{v_1}} \left(\det(\mathbb{M}_{v_1}^1) \dots \det(\mathbb{M}_{v_1}^{d_{v_1}})\right) \prod_{\mathcal{F}_{v_2}} \left(\det(\mathbb{M}_{v_2}^1) \dots \det(\mathbb{M}_{v_2}^{d_{v_2}})\right),$$

where $\mathcal{F}_{v_1}$ ranges over all secondary cells associated to admissible vectors $v_1$, and $\mathcal{F}_{v_2}$ over all cells indexed by non admissible vectors $v_2$.

Choosing all rows and columns of $\mathbb{M}$ indexed by non mixed points, a similar version of Proposition 3.12 holds for the matrix $\mathbb{E}$ instead of $\mathbb{M}$, and it turns out that $\mathbb{E}$ will also have a block structure which will allow us to compute its determinant as follows:

(47)

$$\det(\mathbb{E}) = \det(\mathbb{E}_N) \prod_{\mathcal{F}_{v_1}} \prod_{i=1}^{d_{\tilde{\mathcal{F}}_{v_1}}} \det(\mathbb{E}_{v_1}^i) \prod_{i > d_{\tilde{\mathcal{F}}_{v_1}}} \det(\mathbb{M}_{v_1}^i) \prod_{\mathcal{F}_{v_2}} \prod_i \det(\mathbb{M}_{v_2}^i).$$

From here, the proof of the theorem follows easily, just using Proposition 4.3. □

**Corollary 4.8.** *For every* $i = 0, 1, \dots, n$, *the number of* $i$*-mixed points is exactly* $\mathcal{MV}(Q_0, \dots, Q_{i-1}, Q_{i+1}, \dots, Q_n)$.

**Corollary 4.9.** *If the family* $\{\mathcal{A}_i\}_{i=0,\dots,n}$ *is essential, then* every *coefficient* $c_{i,a}$ *appears in* $\mathrm{Res}_{\mathcal{A}}(f_0, \dots, f_n)$ *with highest power equals to* $\mathcal{MV}(Q_0, \dots, Q_{i-1}, Q_{i+1}, \dots, Q_n)$.

*Remark* 4.10. The same observation given in Remark 3.15 holds also here.

### 4.3. Examples.

**Example 4.11.** Let us compute the resultant given in Example 2.5. In order to do this, we will take $\mathbf{b_0} = (0,0)$ and $Q = \{(0, \frac{1}{3})\}$. Here, $\mathcal{L} = \mathbb{Z}^2$.

The primary cell is $Q_1 + Q_2$, and we will have four secondary cells:

(48)

| $\underline{v}$ | $\underline{type}$ |
|---|---|
| $(2,-1)$ | 2-mixed |
| $(-1,-2)$ | 1-mixed |
| $(-1,-1)$ | non mixed |
| $(-3,-1)$ | 2-mixed. |

In order to subdivide the primary cell, we take $\tilde{\omega}_1 = (0,0,0)$ and $\tilde{\omega}_2 = (1,1)$. This lifting produces two cells: a copy of $Q_1$, and the unique 0-mixed cell of the subdivision. The set $\mathcal{E}$ has 23 points, and we will associate the unique non mixed cell of (48) with $f_2$. With the aid of *Maple*, we have computed $\mathbb{M}$ :

$$\left[\begin{array}{ccccccccccccccccccccccc}
c_1&0&0&0&0&0&c_2&0&0&0&0&0&0&0&0&0&0&0&0&0&0&0&0\\
0&c_1&0&0&0&0&0&0&c_2&0&0&0&0&0&0&0&0&0&0&0&0&0&0\\
0&0&c_1&0&0&0&0&0&0&0&0&0&0&0&0&0&0&c_2&0&0&0&0&0\\
0&0&0&c_1&0&0&0&0&0&0&0&0&0&0&0&0&0&0&c_2&0&0&0&0\\
0&0&0&0&a_1&0&0&0&0&0&0&0&0&0&0&0&a_3&0&0&a_2&0&0&0\\
0&0&0&0&0&a_1&0&0&0&0&0&a_2&0&0&0&0&0&0&0&0&a_3&0&0\\
0&0&0&0&0&0&a_1&0&0&0&0&0&0&0&0&0&0&a_3&0&0&a_2&0&0\\
0&0&a_2&0&0&0&0&a_1&0&0&0&0&0&0&0&0&0&0&0&0&0&a_3&0\\
0&0&0&0&0&0&0&0&a_1&0&0&0&0&0&0&0&0&0&a_3&0&0&a_2&0\\
0&0&0&0&0&0&0&c_2&0&c_1&0&0&0&0&0&0&0&0&0&0&0&0&0\\
0&0&0&0&0&0&0&0&0&0&c_1&0&0&0&0&c_2&0&0&0&0&0&0&0\\
0&0&0&0&0&0&0&0&0&0&0&c_1&0&0&0&0&c_2&0&0&0&0&0&0\\
0&0&0&0&0&0&0&0&0&0&0&0&c_1&0&0&0&0&0&0&c_2&0&0&0\\
0&0&0&0&0&0&0&0&0&0&0&0&0&c_1&0&0&0&0&0&0&c_2&0&0\\
0&0&0&0&0&0&0&0&0&0&0&0&0&0&c_1&0&0&0&0&0&0&c_2&0\\
b_2&0&0&0&b_1&0&0&0&0&0&0&0&0&0&0&b_3&0&0&0&0&0&0&0\\
0&b_2&0&0&0&0&b_1&0&0&0&0&0&0&0&0&0&b_3&0&0&0&0&0&0\\
0&0&0&0&0&0&0&0&b_1&0&0&0&0&0&0&0&0&b_3&0&b_2&0&0&0\\
0&0&0&0&0&0&0&0&0&0&0&0&0&0&0&0&0&0&b_3&0&b_2&0&b_1\\
0&0&0&0&0&b_1&0&0&0&b_2&0&0&0&0&0&0&0&0&0&b_3&0&0&0\\
0&0&0&0&0&0&0&b_1&0&0&b_2&0&0&0&0&0&0&0&0&0&b_3&0&0\\
0&0&0&0&0&0&0&0&0&0&0&b_2&0&0&0&b_1&0&0&0&0&0&b_3&0\\
0&0&0&0&0&0&0&0&0&0&0&0&0&0&0&0&0&0&0&c_1&0&0&c_2
\end{array}\right]$$

Its determinant equals $c_1^4 \operatorname{Res}_{\mathcal{A}}(f_0, f_1, f_2)$. The extraneous factor here is the first principal minor of size $4 \times 4$.

**Example 4.12.** Let us compute the resultant of the same system given in the previous example, but now we will take $Q$ as the unit square $[0,1] \times [0,1]$ translated by the vector $(\epsilon, \frac{1}{3})$, with $\frac{1}{3} \gg \epsilon > 0$. The initial point $\mathbf{b_0}$ will be again $(0,0)$. The primary cell will be now equal to

$Q_1 + Q_2 + Q$, and now we will have 6 secondary cells:

| $v$ | *type* |
|---|---|
| $(2, -1)$ | 2-mixed |
| $(0, -1)$ | non mixed |
| $(-1, -2)$ | 1-mixed |
| $(-1, -1)$ | non mixed |
| $(-3, -1)$ | 2-mixed |
| $(-1, 0)$ | non mixed. |

Observe now that there are secondary cells associated to normal vectors of $Q$ which do not define any facet in $Q_0 + Q_1 + Q_2$. We subdivide the primary cell by taking

$$\begin{array}{llll} \tilde{\omega}_1(0,0) = 1, & \tilde{\omega}_1(1,2) = 1, & \tilde{\omega}_1(2,0) = 0, & \\ \tilde{\omega}_2(1,1) = 1, & \tilde{\omega}_2(3,0) = 0, & & \\ \tilde{\omega}(0,0) = 8, & \tilde{\omega}(1,0) = 4, & \tilde{\omega}(0,1) = 4, & \tilde{\omega}(1,1) = 0. \end{array}$$

If we associate the points lying in the non mixed secondary cells with $f_1$, the matrix $\mathbb{M}$ obtained has size $36 \times 36$ and its determinant equals $c_1^3\, b_3^9\, c_2^3\, b_2^2$ times the resultant. The submatrix $\mathbb{E}$ is the following:

$$\begin{bmatrix}
b_2 & 0 & 0 & b_1 & 0 & 0 & 0 & 0 & 0 & 0 & 0 & 0 & 0 & 0 & 0 & 0 & 0 \\
0 & b_2 & 0 & 0 & 0 & 0 & 0 & 0 & 0 & 0 & 0 & 0 & 0 & 0 & 0 & 0 & 0 \\
0 & 0 & c_1 & 0 & 0 & 0 & 0 & 0 & 0 & 0 & 0 & 0 & 0 & 0 & 0 & 0 & 0 \\
0 & 0 & 0 & c_1 & 0 & 0 & 0 & 0 & c_2 & 0 & 0 & 0 & 0 & 0 & 0 & 0 & 0 \\
0 & 0 & 0 & 0 & c_1 & 0 & 0 & 0 & 0 & 0 & 0 & 0 & 0 & 0 & 0 & 0 & 0 \\
0 & 0 & 0 & 0 & b_2 & b_3 & 0 & 0 & 0 & 0 & 0 & 0 & 0 & 0 & 0 & 0 & 0 \\
0 & 0 & 0 & 0 & 0 & 0 & b_3 & 0 & 0 & 0 & 0 & b_1 & 0 & 0 & 0 & 0 & 0 \\
0 & 0 & b_2 & 0 & 0 & 0 & 0 & b_3 & 0 & 0 & 0 & 0 & 0 & 0 & 0 & 0 & 0 \\
0 & 0 & 0 & 0 & 0 & 0 & 0 & 0 & b_3 & 0 & 0 & 0 & 0 & 0 & 0 & 0 & 0 \\
0 & 0 & 0 & 0 & 0 & 0 & 0 & c_1 & 0 & c_2 & 0 & 0 & 0 & 0 & 0 & 0 & 0 \\
0 & 0 & 0 & 0 & 0 & 0 & 0 & 0 & 0 & 0 & c_2 & 0 & 0 & 0 & 0 & 0 & 0 \\
0 & 0 & 0 & 0 & 0 & 0 & 0 & 0 & 0 & 0 & 0 & c_2 & 0 & 0 & 0 & 0 & 0 \\
0 & 0 & 0 & 0 & 0 & 0 & 0 & 0 & 0 & 0 & 0 & 0 & b_3 & 0 & 0 & 0 & 0 \\
0 & 0 & 0 & 0 & 0 & 0 & 0 & 0 & 0 & b_1 & 0 & 0 & 0 & b_3 & 0 & 0 & 0 \\
0 & 0 & 0 & 0 & 0 & 0 & 0 & b_2 & 0 & 0 & b_1 & 0 & 0 & 0 & b_3 & 0 & 0 \\
0 & 0 & 0 & 0 & 0 & 0 & 0 & 0 & 0 & 0 & 0 & 0 & 0 & 0 & 0 & b_3 & 0 \\
0 & 0 & 0 & 0 & 0 & b_1 & 0 & 0 & 0 & 0 & 0 & 0 & 0 & 0 & 0 & 0 & b_3
\end{bmatrix}$$

**Example 4.13.** We shall see in this example that the inequality established in Proposition 4.3 may be strict. Consider the following essential

family:

$$\begin{array}{rcl} \mathcal{A}_0 & = & \{(0,0);(1,0);(0,1);(1,1)\} \\ \mathcal{A}_1 & = & \{(0,0);(1,1)\} \\ \mathcal{A}_2 & = & \{((1,0);(0,1)\}. \end{array}$$

Denote the generic polynomials having those supports as follows:

$$\begin{array}{rcl} f_0 & = & a_1 + a_2x + a_3y + a_4xy \\ f_1 & = & b_1x + b_2y \\ f_2 & = & c_1 + c_2xy. \end{array}$$

Taking $\mathbf{b_0} = (0,0)$, with $\frac{1}{3} \gg \epsilon > 0$, and $Q$ as the triangle with vertices $(0,0)$, $(1,0)$, $(1,1)$ shifted with $(\epsilon, \frac{1}{3})$, we get

- The primary cell is $Q_1 + Q_2 + Q$.
- There are five secondary cells, associated with the following inward vectors: $(1,-1)$, $(0,-1)$, $(-1,-1)$, $(-1,0)$ and $(-1,1)$.
- $\big(\mathcal{F}_{(1,-1)}\big) \cap \mathbb{Z}^2$ has two points, while an explicit computation reveals that $d_{\tilde{\mathcal{F}}_{(1,-1)}} = 1$.

In this case, $\big(\mathcal{F}_{(1,-1)}\big) \cap \mathbb{Z}^2$ is subdivided in two smaller cells, each of them contains exactly one point, and any of the two points may be chosen as mixed. If we take $\tilde{\omega}_i$, $i = 1,2$ always equal to 2, $\tilde{\omega}$ always equal to 0, and apply the algorithm given in the paragraph 4.1, we will get the following $13 \times 13$ matrix:

$$\mathbb{M} = \begin{bmatrix} c_2 & 0 & 0 & 0 & 0 & 0 & 0 & 0 & 0 & 0 & 0 & c_1 & 0 \\ 0 & c_2 & 0 & 0 & 0 & 0 & 0 & 0 & 0 & 0 & 0 & 0 & c_1 \\ b_2 & 0 & b_1 & 0 & 0 & 0 & 0 & 0 & 0 & 0 & 0 & 0 & 0 \\ 0 & 0 & 0 & b_1 & 0 & 0 & 0 & 0 & 0 & 0 & 0 & b_2 & 0 \\ 0 & 0 & 0 & 0 & b_1 & 0 & 0 & 0 & b_2 & 0 & 0 & 0 & 0 \\ 0 & 0 & 0 & 0 & 0 & b_1 & 0 & 0 & 0 & b_2 & 0 & 0 & 0 \\ 0 & b_1 & 0 & 0 & 0 & 0 & b_2 & 0 & 0 & 0 & 0 & 0 & 0 \\ 0 & 0 & 0 & 0 & 0 & 0 & 0 & b_2 & 0 & b_1 & 0 & 0 & 0 \\ c_1 & 0 & 0 & 0 & 0 & 0 & 0 & 0 & c_2 & 0 & 0 & 0 & 0 \\ 0 & c_1 & 0 & 0 & 0 & 0 & 0 & 0 & 0 & c_2 & 0 & 0 & 0 \\ 0 & b_2 & 0 & 0 & 0 & 0 & 0 & 0 & 0 & 0 & b_1 & 0 & 0 \\ a_4 & a_2 & 0 & 0 & 0 & 0 & a_3 & 0 & 0 & 0 & 0 & a_1 & 0 \\ 0 & a_4 & 0 & a_2 & 0 & 0 & 0 & 0 & 0 & 0 & 0 & a_3 & a_1 \end{bmatrix}.$$

We have ordered the rows and columns of $\mathbb{M}$ in such a way that the monomials lying in $\big(\mathcal{F}_{(1,-1)}\big) \cap \mathbb{Z}^2$ index rows and columns 7 and 8. All the other non mixed points index the first six rows and columns of $\mathbb{M}$. An explicit computation reveals that

$$\det(\mathbb{M}) = b_2\, c_2^2\, b_1^4\, \mathrm{Res}_{\mathcal{A}}(f_0, f_1, f_2)$$

and it is straightforward to check that, in this case, $b_2c_2^2b_1^4$ may be obtained by taking the principal minor of size $7 \times 7$, or computing the minor indexed by the monomials lying in the first six rows and the monomial number eight.

Departamento de Matemática, F.C.E.y N., Universidad de Buenos Aires. Pabellón I, Ciudad Universitaria (1428) Buenos Aires, Argentina.
*E-mail address*: cdandrea@dm.uba.ar